\documentclass[twoside]{article}

\usepackage[utf8]{inputenc}
\usepackage{graphicx}
\usepackage[all]{xy}
\usepackage{amssymb}
\usepackage{amsmath} 
\usepackage{amsthm}
\usepackage{verbatim}
\usepackage{latexsym}
\usepackage{mathrsfs}
\usepackage{makeidx}
\usepackage{setspace}

\usepackage{mathrsfs}
\usepackage{watermark}

\newenvironment{myItemize}{ 
  \begin{list}{\raisebox{2.2pt}{$\centerdot$}}{%
      \setlength\leftmargin{18pt}
      \setlength\labelwidth{20pt}
    }
  }{
  \end{list}
}

\newenvironment{myEnumerate}{ 
  \begin{list}{\labelenumi}{%
      \usecounter{enumi}
      \setlength\leftmargin{25pt}
      \setlength\labelwidth{20pt}
    }
  }{
  \end{list}
}

\newenvironment{myAlphanumerate}{ 
  \begin{list}{\labelenumii}{%
      \usecounter{enumii}
      \setlength\leftmargin{25pt}
      \setlength\labelwidth{20pt}
    }
  }{
  \end{list}
}

\usepackage[T1]{fontenc}

\usepackage{epigraph}

\renewcommand{\Cup}{\bigcup}
\renewcommand{\Cap}{\bigcap}

\renewcommand{\a}{\alpha}
\renewcommand{\b}{\beta}
\newcommand{\g}{\gamma}

\renewcommand{\d}{\delta}
\renewcommand{\k}{\kappa}
\renewcommand{\l}{\lambda}

\renewcommand{\o}{\omega}
\newcommand{\s}{\sigma}

\newcommand{\nr}{{\operatorname{nr}}}

\newcommand{\es}{\varnothing}

\newcommand{\B}{\mathcal{B}}

\newcommand{\D}{\mathcal{D}}
\newcommand{\E}{\mathcal{E}}

\newcommand{\R}{\mathbb{R}} 

\newcommand{\cat}{{}^\frown}  

\newcommand{\la}{\langle} 
\newcommand{\ra}{\rangle}
\newcommand{\Po}{\mathcal{P}}             
 
\newcommand{\NS}{{\operatorname{NS}}}

\newcommand{\cf}{\operatorname{cf}}

\newcommand{\reg}{\operatorname{reg}} 
\newcommand{\id}{\operatorname{id}}

\newcommand{\GC}{\operatorname{GC}}

\newcommand{\rest}{\!\restriction\!}
\newcommand{\restl}{\restriction}  
\newcommand{\dom}{\operatorname{dom}}

\newcommand{\OTP}{\operatorname{OTP}}
  
\renewcommand{\le}{\leqslant}  
\renewcommand{\ge}{\geqslant}  
\newcommand{\less}{\lessdot}  
\newcommand{\sd}{\,\triangle\,}             

\newcommand{\Pii}{{\Pi_1^1}}

\newcommand{\PlOne}{\,{\textrm{\bf I}}}
\newcommand{\PlTwo}{\textrm{\bf I\hspace{-1pt}I}}

\swapnumbers
\newtheorem*{Thm*}{Theorem}
\newtheorem{Thm}{Theorem}
\newtheorem{Lemma}[Thm]{Lemma}
\newtheorem{Cor}[Thm]{Corollary}

\newenvironment{claim}[1]{\text{ }\vspace{7pt}\newline\noindent\textbf{Claim #1.}}{\hspace*{\fill}}
\newenvironment{claim*}{\vspace{7pt}\noindent\textbf{Claim.}}{}

\theoremstyle{definition}
\newtheorem{Def}[Thm]{Definition}
\newtheorem*{Def*}{Definition}

\newtheorem*{Example*}{Example}
\newtheorem*{Notation}{Notation}
\newtheorem*{Assumption}{Assumption}

\theoremstyle{remark}
\newtheorem*{Remark}{Remark}
\newtheorem{RemarkN}[Thm]{Remark}

\author{Vadim Kulikov\thanks{EDIT 2025: The name of the author has changed to Vadim Weinstein} \\
  Department of Mathematics and Statistics\\ University of Helsinki}
\title{Borel Reductions and Cub Games in Generalized Descriptive Set Theory}

\date{July 12, 2012\\(file edited \today)}

\begin{document}

\maketitle

 \begin{abstract}
   It is shown that the power set of $\k$ ordered by the subset relation modulo various versions of the
   non-stationary ideal can be embedded into the partial order of Borel equivalence relations
   on $2^\k$ under Borel reducibility. Here $\k$ is uncountable regular cardinal
   with $\k^{<\k}=\k$.
 \end{abstract}

\subsection*{Acknowledments}

This work is part of my Ph.D. thesis, to be defended in fall 2011. 
I wish to express my gratitude to my supervisor Tapani Hyttinen
for his careful attention towards this work, his valuable advice
and some technical help regarding the 
proof of Theorem~\ref{thm:Main2}.

I am indebted to
the Finnish National Graduate School in Mathematics and its Applications
for supporting my post-graduate studies during the preparation of this article.

\section{Introduction}

It is shown that the partial order of Borel
equivalence relations on the generalized Baire spaces ($2^\k$ for $\k>\o$) 
under Borel reducibility has high complexity already at low levels (below $E_0$).

This extends an answer stated in~\cite{5HytFri} to an open problem stated in~\cite{5FriHytKul} and in particular
solves open problems 7 and 9 from~\cite{5HytFri}.

The developement of the theory 
of the generalized Baire and Cantor spaces dates back to 1990's when it
A. Mekler and J. Väänänen published the paper \emph{Trees and $\Pii$-Subsets of ${}^{\o_1}\o_1$} \cite{4MekVaa}
and A. Halko published \emph{Negligible subsets of the generalized Baire space $\o_1^{\o_1}$}. 
More recently equivalence relations and Borel reducibility 
on these spaces and their applications to model theory 
have been under focus, see my latest joint work 
with S. Friedman and T. Hyttinen~\cite{5FriHytKul}.

Suppose $\k$ is an infinite cardinal and let $\E^B_{\k}$ be the collection of all 
Borel equivalence relations
on $2^\k$. (For definitions in the case $\k>\o$ see next section.) 
For equivalence relations $E_0$ and $E_1$ let us denote $E_0\le_B E_1$  if there exists
a Borel function $f\colon 2^\k\to 2^\k$ such that $(\eta,\xi)\in E_0\iff (f(\eta),f(\xi))\in E_1$.
The relation $\le_B$ defines a quasiorder on $\E^B_\k$, i.e. it induces a partial order
on $\E^B_\k/\sim_B$ where $\sim_B$ is the equivalence relation of bireducibility: 
$E_0\sim_B E_1\iff (E_0\le_B E_1)\land (E_1\le_B E_0)$.

In the case $\k=\o$ there are many known results 
that describe the order $\la\E^B_\k,\le_B\ra$. 
Two of them are:

\begin{Thm*}[Louveau-Velickovic \cite{4LouVel}]
  The partial order $\la \Po(\o),\subset_*\ra$ can be embedded into
  the partial order \mbox{$\la \E^{B}_\o,\le_B\ra$}, where 
  $A\subset_* B$ if $A\setminus B$ is finite.
\end{Thm*}

\begin{Thm*}[Adams-Kechris \cite{4AdaKec}]
  The partial order $\la \B,\subset\ra$ can be embedded into the partial order 
  \mbox{$\la \E^{B}_\o,\le_B\ra$}, where 
  $\B$ is the collection of all Borel subsets of the real line $\R$. In fact, the embedding is into 
  the suborder of $\la \E^{B}_\o,\le_B\ra$
  consisting of the countable Borel equivalence relations, i.e., those Borel equivalence relations each 
  of whose equivalence classes is countable.
\end{Thm*}

Our aim is to generalize these results to uncountable $\k$ with $\k^{<\k}=\k$
and it is proved that $\la \Po(\k),\subset_{\NS(\o)}\ra$ can be embedded into 
$\la \E^{B}_\k,\le_B\ra$, where  $A\subset_{\NS(\o)} B$ means that $A\setminus B$ is not $\o$-stationary.
This is proved in ZFC. However under mild additional assumptions on $\k$ or on the underlying set theory,
it is shown that $\la \Po(\k),\subset_\NS\ra$ can be embedded into 
$\la \E^{B}_\k,\le_B\ra$, where  $A\subset_\NS B$ means that $A\setminus B$ is non-stationary
and that $\la \Po(\k),\subset_*\ra$ can be embedded into 
$\la \E^{B}_\k,\le_B\ra$, where  $A\subset_* B$ means that $A\setminus B$ is bounded.

\begin{Assumption}
  Everywhere in this article it is assumed that $\k$ is a cardinal which satisfies 
  $|\k^{\a}|=\k$ for all $\a<\k$. This requirement is briefly denoted by $\k^{<\k}=\k$.
\end{Assumption}

\section{Background in Generalized Descriptive Set~Theory}

\begin{Def}
  Consider the function space $2^\k$ (all functions from $\k$ to $\{0,1\}$) 
  equipped with the
  topology generated by the sets 
  $$N_p=\{\eta\in 2^\k\mid \eta\rest \a=p\}$$ 
  for $\a<\k$ and $p\in 2^{\a}$.
  Borel sets on this space are obtained by closing the topology under unions and intersections of 
  length~$\le\k$, and complements. 

  An equivalence relation $E$ on $2^\k$ is \emph{Borel reducible} 
  to an equivalence relation $E'$  on $2^\k$ if
  there exists a Borel function $f\colon 2^\k\to 2^\k$ (inverse images of open sets are Borel)
  such that $\eta E\xi\iff f(\eta)E' f(\xi)$. This
  is denoted by $E\le_B E'$.

  The descriptive set theory of these spaces, of
  equivalence relations on them and of their reducibility properties for $\k>\o$, has been developed at least 
  in~\cite{5FriHytKul,4Halko,4MekVaa}. For $\k=\o$ this is the field of standard descriptive set theory.
  
  By $\id_{X}$ we denote the identity relation on $X$:
  $(\eta,\xi)\in \id_{X}\iff (\eta,\xi)\in X^2\land \eta=\xi$ and by $E_0$ the equivalence relation on $2^\k$ 
  (or on $\k^\k$ as in the proof of Theorem~\ref{thm:E_SEmbedstoE_0})
  such that $(\eta,\xi)\in E_0\iff \{\a\mid \eta(\a)\ne\xi(\a)\}\text{ is bounded}$.
\end{Def}

\begin{Notation}
  Let $\E^B_\k$ denote the set of all Borel equivalence relations on $2^\k$
  (i.e. equivalence relations $E\subset (2^\k)^2$ such that $E$ is 
  a Borel set). If $X,Y\subset \k$ and $X\setminus Y$ is non-stationary, let us
  denote it by $X\subset_{\NS} Y$. If $X\setminus Y$ is not $\l$-stationary for
  some regular $\l<\k$, it is denoted by $X\subset_{\NS(\l)} Y$.

  The set of all ordinals below $\k$ which have cofinality $\l$ is denoted by $S^\k_\l$,
  and $\lim (\k)$ denotes the set of all limit ordinals below~$\k$. Also $\reg\k$ denotes the set of
  regular cardinals below $\k$ and 
  $$S^\k_{\ge\l}=\Cup_{{\mu\ge\l}\atop{\mu\in\reg\k}}S_{\mu}^\k,$$
  $$S^\k_{\le\l}=\Cup_{{\mu\le\l}\atop{\mu\in\reg\k}}S_{\mu}^\k.$$

  If $A\subset \a$ and $\a$ is an ordinal, then $\OTP(A)$ is the order type of $A$ in the ordering induced on it 
  by~$\a$.\label{page:OTP}

  For ordinals $\a<\b$ let us adopt the following abbreviations:
  \begin{myItemize}
  \item $(\a,\b)=\{\g\mid \a<\g<\b\}$,
  \item $[\a,\b]=\{\g\mid \a\le\g\le\b\}$,
  \item $(\a,\b]=\{\g\mid \a<\g\le\b\}$,
  \item $[\a,\b)=\{\g\mid \a\le\g<\b\}$.
  \end{myItemize}

  If $\eta$ and $\xi$ are functions in $2^\k$, then $\eta\sd\xi$ is the function
  $\zeta\in 2^\k$ such that $\zeta(\a)=1\iff \eta(\a)\ne\xi(\a)$ for all $\a<\k$, and $\bar\eta=1-\eta$ is the function $\zeta\in 2^\k$
  such that $\zeta(\a)=1-\eta(\a)$ for all $\a<\k$. If $A$ and $B$ are sets,
  then $A\sd B$ is just the symmetric difference.

  For any set $X$, $2^{X}$ denotes the set of all functions from $X$ to $2=\{0,1\}$.
  If $p\in 2^{[0,\a)}$ and $\eta\in 2^{[\a,\k)}$, then $p\cat \eta\in 2^\k$ is the catenation:
  $(p\cat\eta)(\b)=p(\b)$ for $\b<\a$ and $(p\cat\eta)(\b)=\eta(\b)$ for $\b\ge\a$.
\end{Notation}

\begin{Def}\label{def:CoMeagerStuff}
  A \emph{co-meager} subset of $X$ is a set 
  which contains an intersection of length $\le\k$
  of dense open subsets of $X$. Co-meager sets are always non-empty and form a filter on $2^\k$,~\cite{4MekVaa}.
  A set $X$ has the \emph{Property of Baire} if there exists an open set $A$ such that $X\sd A$
  is meager, i.e. a complement of a co-meager set. 
  As in standard descriptive set theory, Borel sets have the Property of Baire (proved in \cite{4Halko}). 
  For a Borel function $f\colon 2^\k\to 2^\k$ denote by $C(f)$ one of the co-meager sets 
  restricted to which $f$ is continuous
  (such set is not unique, but we can always pick one using the Property of Baire of Borel sets, 
  see \cite{5FriHytKul}).
\end{Def}

\begin{Lemma}\label{lem:YouCanFind}
  Let $D$ be a co-meager set in $2^\k$ and let $p,q\in 2^\a$ for some $\a<\k$. Then
  there exists $\eta\in 2^{[\a,\k)}$ such that $p\cat\eta\in D$ and $q\cat\eta\in D$.
  Also there exists $\eta\in 2^{[\a,\k)}$ such that $p\cat\bar\eta\in D$ and $q\cat\eta\in D$
  where $\bar\eta=1-\eta$.
\end{Lemma}
\begin{proof}
  Let $h$ be the homeomorphism $N_p\to N_q$ defined by $p\cat \eta\mapsto q\cat\eta$. 
  Then $h[N_p\cap D]$ is co-meager in
  $N_q$, so $N_q\cap D\cap h[N_p\cap D]$ is non-empty. Pick $\eta'$ from that intersection
  and let $\eta=\eta'\rest[\a,\k)$. This will do. For the second part take for $h$ the homeomorphism
  defined by $p\cat\eta\mapsto q\cat\bar\eta$.
\end{proof}

\section{On Cub-games and $\GC_\l$-characterization}\label{sec:GClchar}

The notion of cub-games is a useful way to treat certain properties of subsets of cardinals. 
They generalize closed unbounded sets and are related to combinatorial principles 
such as $\square_{\k}$. Under mild set theoretic assumptions, they give characterizations of
CUB-filters in different cofinalities. Treatments of this subject can be found 
for example in~\cite{5HuuHytRau, 3Hyttinen, 3HytSheTuu}.

\begin{Def}\label{def:cubgame2}
  Let $A\subset \k$. The game $\GC_{\l}(A)$ is played between players $\PlOne$
  and $\PlTwo$ as follows. There are $\l$ moves and at the $i$:th 
  move player $\PlOne$ picks an ordinal $\a_i$ which is greater than any ordinal
  picked earlier in the game and then $\PlTwo$ picks an ordinal $\b_i>\a_i$.
  Player $\PlTwo$ wins if $\sup_{i<\l}\a_i\in A$. Otherwise player $\PlOne$ wins.  
\end{Def}

\begin{Def}
  A set $C\subset \k$ is \emph{$\l$-closed} for a regular cardinal $\l<\k$, if 
  for all increasing sequences $\la\a_i\in C\mid i<\l\ra$, the limit $\sup_{i<\l}\a_i$ is in $C$.
  A set $C\subset \k$ is \emph{closed} if it is $\l$-closed for all regular $\l<\k$.
  A set is \emph{$\l$-cub} if it is $\l$-closed and unbounded and \emph{cub}, if it is
  closed and unbounded. A set is $\l$-stationary, if it intersects all $\l$-cub sets and
  \emph{stationary} if it intersects all cub sets.
\end{Def}

\begin{Def}
  We say that $GC_\l$-characterization holds for $\k$, if 
  $$\{A\subset \k\mid \PlTwo\text{ has a winning strategy in }\GC_\l(A)\}
  =\{A\subset\k\mid A\text{ contains a }\l\text{-cub set}\}$$
  and we say that $\GC$-characterization holds for $\k$ if 
  $\GC_\l$-characterization holds for $\k$ for all regular $\l<\k$.
\end{Def}

\begin{Def}
  Assume $\k=\l^+$ and $\mu\le\l$ a regular uncountable cardinal.  The
  \emph{square principle on $\k$ for $\mu$}, denoted $\square^\k_\mu$,
  defined by Jensen in case $\l=\mu$, is the statement that there
  exists a sequence $\la C_\a\mid\a\in S^\k_{\le\mu}\ra$ with the
  following properties:
  \begin{myEnumerate}
  \item $C_\a\subset\a$ is closed and unbounded in $\a$,
  \item if $\b\in \lim C_\a$, then $C_\b=\b\cap C_\a$,
  \item if $\cf(\a)<\mu$, then $|C_\a|<\mu$.
  \end{myEnumerate}
\end{Def}
\begin{RemarkN}\label{rem:Shelahonsquares}
  For $\o<\mu<\l$ in the definition above, it was proved by Shelah in \cite{5Shelah2} that $\square^{\k}_\mu$ holds
  (can be proved in ZFC, for a proof see also \cite[Lemma 7.7]{3BurkeMagidor}).
  If $\mu=\l$, then $\square^\k_\mu=\square_\mu^{\mu^+}$ is denoted by $\square_\mu$ and
  can be easily forced or, on the other hand, it holds, if $V=L$.
  The failure of $\square_\mu$ implies that $\mu^+$ is Mahlo in $L$, as pointed 
  out by Jensen, see~\cite{4Jech}. 
\end{RemarkN}

\begin{Def}
  For $\k>\o$, the set $I[\k]$ 
  consists of those $S\subset \k$ that have the following property:
  there exists a cub set $C$ and a sequence $\la \D_\a\mid \a<\k\ra$ such that
  \begin{myEnumerate}
  \item $\D_\a\subset \Po(\a)$, $|\D_\a|<\k$,
  \item $\D_\a\subset \D_\b$ for all $\a<\b$,
  \item for all $\a\in C\cap S$ there exists $E\subset \a$ unbounded in $\a$ and of order type $\cf(\a)$
    such that for all $\b<\a$, $E\cap \b\in \D_\g$ for some $\g<\a$.
  \end{myEnumerate}
\end{Def}

\begin{RemarkN}\label{rem:KnownGC}
  The following is known.
  \begin{myEnumerate}
  \item $I[\k]$ is a normal ideal and contains the non-stationary sets.
  \item If $\l<\k$ is regular and $S^\k_\l\in I[\k]$, then $\GC_\l$-characterization holds for $\k$.
  \item If $\mu$ is regular and $\k=\mu^+$, then $S^\k_{<\mu}\in I[\k]$,~\cite{5Shelah2}.
    This follows also from 4. and Remark \ref{rem:Shelahonsquares}
  \item When $\l>\o$, then $\square^\k_\l$ implies that $S^\k_\l\in I[\k]$ (take $\D_\a=\{C_\a\cap\b\mid\b<\a\}$).
  \item $S^\k_\o\in I[\k]$.
  \item If $\k^{<\l}=\k=\l^+$, then $\GC_\l$-characterization holds for $\k$ if and only if
    $\k\in I[\k]$ if and only if $S^\k_\l\in I[\k]$, see~\cite[Corollary~2.4]{5HuuHytRau} and~\cite{5Shelah2}.
  \item The existence of $\l<\k$ such that
    $\GC_\l$-characterization does not hold for $\k$ is equiconsistent with the 
    existence of a Mahlo 
    cardinal.\footnote{A good exposition of this result can be found in Lauri Tuomi's Master's thesis 
      (University of Helsinki, 2009).}
    Briefly this is because the failure of the characterization implies the failure of $\square_{\l}$
    which implies that $\l^+$ is Mahlo in $L$ as discussed above.
    On the other hand,
    in the Mitchell model, obtained from $S_{\text{in}}=\{\d<\l\mid \d\text{ is inaccessible}\}$ where $\l>\k$ is Mahlo,
    it holds that $S_{\text{in}}\notin I[\k^+]$, \cite[Lemma~2.6]{5HuuHytRau}.

  \item \label{remit5} If $\k$ is regular and for all regular $\mu<\k$ we have $\mu^{<\l}<\k$, 
    then $\k\in I[\k]$.
  \end{myEnumerate}
\end{RemarkN}

\begin{Remark}
  As Remark \ref{rem:KnownGC} shows, the assumption that $\GC_\l$-characterization holds for $\k$ is quite weak.
  For instance $\GC_\o$-characterization holds for all regular $\k>\o$ and GCH implies that
  $\GC_\l$-characterization holds for $\k$
  for all regular $\l<\k$. 
\end{Remark}

\section{Main Results}

Theorems \ref{thm:Main1} and \ref{thm:Main2} constitute the goal of this work.
They are stated below but proved in the end of this section, starting at
pages \pageref{pageproof1} and \pageref{pageproof2} respectively.

\begin{Thm}\label{thm:Main1}
  Assume that $\l<\k$ are regular and $\GC_\l$-characterization holds for $\k$. Then 
  the order $\la\Po(\k),\subset_{\NS(\l)}\ra$ can be embedded into $\la \E^B_\k,\le_B\ra$ 
  strictly between $\id_{2^\k}$ and $E_0$. More precisely
  there exists a one-to-one map $F\colon \Po(\k)\to \E^{B}_\k$ such that
  for all $X,Y\in \Po(\k)$ we have $\id_{2^\k}\lneq_B F(X)\lneq_B E_0$ and 
  $$X\subset_{\NS(\l)} Y\iff F(X)\le_B F(Y).$$
\end{Thm}

\begin{Thm}\label{thm:Main2}
  Assume either $\k=\o_1$ or $\k=\l^+>\o_1$ and $\ \square_\l$. Then the
  partial order \mbox{$\la \Po(\k),\subset_{\NS}\ra$} can be embedded into $\la \E^B_\k,\le_B\ra$.
\end{Thm}

\subsection{Corollaries}

\begin{Cor}
  Assume that $\l<\k$ is regular. Additionally assume one of the following:
  \begin{myEnumerate}
  \item $\k=\mu^+$, $\mu$ is regular and $\l<\mu$,
  \item $\k=\l^+$ and $\square_\l$ holds,
  \item for all regular $\mu<\k$, $\mu^{<\l}<\k$ (e.g. $\k$ is $\o_1$ or inaccessible).
  \end{myEnumerate}
  Then the
  partial order $\la \Po(\k),\subset_{\NS(\l)}\ra$ can be embedded into $\la \E^B_\k,\le_B\ra$.
\end{Cor}
\begin{proof}
  Any of the assumptions 1 -- 4 is sufficient to obtain $\GC_\l$-characterization for $\k$
  by Remarks \ref{rem:KnownGC} and~\ref{rem:Shelahonsquares}, so 
  the result follows from Theorem~\ref{thm:Main1}.
\end{proof}

\begin{Cor}\label{cor:First}
  The partial order $\la \Po(\k),\subset_{\NS(\o)}\ra$ can be embedded into \mbox{$\la \E^B_\k,\le_B\ra$}.
  In particular $\la \Po(\o_1),\subset_{\NS}\ra$ can be embedded into $\la \E^B_{\o_1},\le_B\ra$ 
  assuming CH.
\end{Cor}
\begin{proof}
  By Remark \ref{rem:KnownGC} $\GC_\o$-characterization holds for 
  $\k$ for any regular $\k>\o$, so the result follows from Theorem~\ref{thm:Main1}.
\end{proof}

\begin{Def}
  Let $S\subset\k$. Then the combinatorial principle $\diamondsuit_\k(S)$ 
  states that
  there exists a sequence $\la D_\a\mid \a\in S\ra$ such that for every $A\subset \k$
  the set $\{\a\mid A\cap \a=D_\a\}$ is stationary.
\end{Def}

\begin{Thm}[Shelah \cite{4Shelah6}]\label{thm:diamond}
  If $\k=\l^+=2^\l$
  and $S\subset \k\setminus S^\k_{\cf(\l)}$ is stationary, then
  $\diamondsuit_\k(S)$ holds. \qed
\end{Thm}

\begin{Cor}\label{cor:Second}
  \begin{myEnumerate}
  \item The ordering $\la\Po(\k),\subset\ra$ can be embedded into $\la \E^{B}_{\k},\le_B\ra$.
  \item Assume that $\k=\o_1$ and $\diamondsuit_{\o_1}$ holds or that 
    $\k$ is not a successor of an $\o$-cofinal cardinal.
    Then also the ordering $\la \Po(\k),\subset_*\ra$ can be embedded
    into \mbox{$\la \E^{B}_{\k},\le_B\ra$}, where $\subset_*$ is 
    inclusion modulo bounded sets.
  \end{myEnumerate}
\end{Cor}
\begin{proof}
  For the first part it is sufficient to show that 
  the partial order \mbox{$\la\Po(\k),\subset\ra$} can be embedded
  into \mbox{$\la\Po(\k),\subset_{\NS(\o)}\ra$}. 
  Let $G(A)=\Cup_{i\in A}S_i$ where $\{S_i\subset S^\k_\o\mid i<\k\}$
  is a collection of disjoint stationary sets. Then $A\subset B\iff G(A)\subset_{\NS} G(B)$,
  so this proves the first part.

  For the second part, let us show that if $\diamondsuit_\k(S^\k_\l)$ holds, then
  $\la \Po(\k),\subset_*\ra$ can be embedded into \mbox{$\la\Po(\k),\subset_{\NS(\l)}\ra$}.
  Then the result follows. If $\k=\o_1$ and $\diamondsuit_{\o_1}$ holds, then it follows
  by Corollary~\ref{cor:First}. On the other hand, if $\k$ is not a successor of an $\o$-cofinal cardinal, then
  from Theorem~\ref{thm:diamond} it follows that $\diamondsuit_\k(S^{\k}_\o)$ holds
  and the result follows again from Corollary~\ref{cor:First}.

  Suppose that $\la D_\a\mid \a\in S^\k_\l\ra$ is a $\diamondsuit_\k(S^\k_\l)$-sequence. 
  If $X,Y\subset \a$ for $\a\le \k$, let $X\subset_* Y$ denote that there is $\b<\a$
  such that $X\setminus \b\subset Y\setminus \b$, i.e. $X$ is a subset of $Y$ on a final segment of $\a$.
  Note that this coincides with the earlier defined $\subset_*$ when $\a=\k$.
  For $A\subset \k$
  let
  $$H(A)=\{\a<\k\mid D_\a\subset_* A\cap \a\}.$$
  If $A\subset_*B$ then there is $\g<\k$ such that $A\setminus \g\subset B\setminus \g$
  and if $\b>\g$ is in $H(A)$, then $D_\b\subset_* A\cap \b$ and since $A\cap\b\subset_* B\cap\b$,
  we have $D_\b\subset_* B\cap\b$, so $H(A)\subset_* H(B)$ which finally implies $H(A)\subset_{\NS(\o)}H(B)$.

  Assume now that $A\not\subset_* B$ and let $C=A\setminus B$. Let $S'$ be the stationary
  set such that for all $\a\in S'$, $C\cap \a=D_\a$. Let $S$ be the $\l$-stationary set
  $S'\cap \{\a\mid C\text{ is unbounded below }\a\}$. $S$ is stationary, because it 
  the intersection of $S'$ and a cub set. Now for all $\a\in S$ we have  
  $D_\a=C\cap\a\subset A\cap\a$, so $S\subset H(A)$. On the other hand if $\a\in S$, then
  $$D_\a\setminus (B\cap \a)=(C\cap\a)\setminus (B\cap \a)=
    ((A\setminus B)\cap\a)\setminus (B\cap \a)=C\cap \a$$
  is unbounded in $\a$, so $D_\a\not\subset_* B\cap\a$ and so $S\subset H(A)\setminus H(B)$,
  whence $H(A)\not\subset_{\NS(\l)} H(B)$.
\end{proof}

\begin{Cor}
  There are $2^\k$ equivalence relations between $\id$ and $E_0$ that form a linear order
  with respect to $\lneq_B$.
\end{Cor}
\begin{proof}
  Let $K=\{\eta\in 2^\k\mid (\exists\b)(\forall\g>\b)(\eta(\g)=0)\}$,
  let $f\colon K \to \k$ be a bijection and for $\eta,\xi\in 2^\k$ define
  $\eta\less \xi$ if and only if 
  $$\eta(\min\{\a\mid \eta(\a)\ne \xi(\a)\})<\xi(\min\{\a\mid \eta(\a)\ne \xi(\a)\}).$$
  For $\eta\in 2^\k$ let $A_\eta=\{f(\xi)\mid \xi\less \eta\land \xi\in K\}$.
  Clearly $A_{\eta}\subsetneq A_\xi$ if and only if $\eta\less \xi$ and the latter is a linear order.
  The statement now follows from Corollary~\ref{cor:Second}.
\end{proof}

\subsection{Preparing for the Proofs}

\begin{Def}\label{def:EqRels}
  For each $S\subset \lim\k$ 
  let us define equivalence relations $E_S^*$,  $E_S$
  and $E_S^*(\a)$,  $\a\le\k$,
  on the space $2^{\k}$ as follows. 
  Suppose $\eta,\xi\in 2^\d$ for some $\d\le\k$ and let
  $\zeta=\eta\sd \xi$. Let us 
  define $\eta$ and $\xi$ to be $E_S^*(\d)$-equivalent if and only if for all ordinals $\a\in S\cap\d$ 
  there exists $\beta<\a$ such that $\zeta(\g)$ has the same value for all $\g\in(\b,\a)$.
  Let $E_S^*=E_S^*(\k)$ and $E_S=E_S^*\cap E_0$, where $E_0$ is the equivalence modulo
  bounded sets.
\end{Def}

\begin{Remark}
  If $S=\es$, then $E_S=E_\es=E_0$. If $S=\lim\k$ or equivalently if $S=\lim_{\o}\k=S^\k_\o$ 
  ($\o$-cofinal limit ordinals),
  then $E_S=E_0'$, where $E_0'$ is defined in~\cite{5HytFri}.
\end{Remark}

\begin{Thm}
  For any $S\subset \lim\k$ the equivalence relations $E_S$ and $E_S^*$ are Borel.
\end{Thm}
\begin{proof} This is obvious by writing out the definitions:
  \begin{eqnarray*} 
    E_S^*&=&\Cap_{\a\in S}\Cup_{\b<\a} \Big(\Cap_{\b<\g<\a}\{(\eta,\xi)\mid \eta(\g)\ne\xi(\g)\}\cup\Cap_{\b<\g<\a}\{(\eta,\xi)\mid \eta(\g)=\xi(\g)\}\Big),\\
    E_0&=&\Cup_{\a<\k}\Cap_{\a<\b<\k}\{(\eta,\xi)\mid \eta(\b)=\xi(\b)\}.\\
    E_S&=&E_S^*\cap E_0.
  \end{eqnarray*}
\end{proof}

The ideas of the following proofs are simple, but are repeated many times in this article
in one way or another.

\begin{Thm}\label{thm:NotToID}
   For all $S\subset \lim\k$, $E_S\not\le_B\id_{2^\k}$ and $E_S^*\le_B\id_{2^\k}$.
\end{Thm}
\begin{proof}
  For the first part suppose $f$ is a Borel reduction from $E_S$ to $\id_{2^\k}$. 
  Let $\eta$ be a function such that $\eta$ and $\bar\eta=1-\eta$
  are both in $C(f)$ (see Definition \ref{def:CoMeagerStuff}, page~\pageref{def:CoMeagerStuff}).
  This is possible by Lemma \ref{lem:YouCanFind}, page~\pageref{lem:YouCanFind}.
  Then $(\eta,\bar\eta)\notin E_S$. Let $\a$ be so large that
  $f(\eta)\rest \a\ne f(\bar\eta)\rest\a$ and pick $\b$ so that
  $$f[N_{\eta\restl\b}\cap C(f)]\subset N_{f(\eta)\restl\a}$$
  and
  $$f[N_{\bar\eta\restl\b}\cap C(f)]\subset N_{f(\eta)\restl\a}.$$
  This is possible by the continuity of $f$ on $C(f)$.
  By Lemma \ref{lem:YouCanFind} pick now a $\zeta\in 2^{[\beta,\k)}$
  so that $\eta\rest\b\cat\zeta\in C(f)$ and $\bar\eta\rest\b\cat\zeta\in C(f)$
  which provides us with a contradiction, since
  $$\big(\eta\rest\b\cat\zeta, \bar\eta\rest\b\cat\zeta\big)\in E_S,
     \text{ but }f(\eta\rest\b\cat\zeta)\ne f(\bar\eta\rest\b\cat\zeta)$$

  To prove the second part it is sufficient to construct a reduction from $E_S^*$ to $\id_{\k^\k}$,
  since $\id_{\k^\k}$ and $\id_{2^\k}$ are bireducible (see \cite{5FriHytKul}).
  Let us define an equivalence relation $\sim$ on $2^{<\k}$ such that $p\sim q$ if and only if 
  $\dom p=\dom q$ and $p\sd q$ is eventually
  constant, i.e. for some $\a<\dom p$, $(p\sd q)(\g)$ is the same for all $\g\in [\a,\dom p)$.\label{pagesimrel}
  Let $s\colon 2^{<\k}\to \k$ be a map such that $p\sim q\iff s(p)=s(q)$. Suppose $\eta\in 2^\k$ and let us define
  $\xi=f(\eta)$ as follows. Let $\b_\g$ denote the $\g$:th element of $S$ and
  let $\xi(\g)=s(\eta\rest\b_\g)$. Now we have $\eta E_S^* \xi$ if and only if $\eta\rest\b_\g=\xi\rest\b_\g$
  for all $\g\in \k$ if and only if $f(\eta)=f(\xi)$. 
\end{proof}

\begin{Cor}\label{cor:NotToID2}
  Let $S\subset \k$. 
  If $p\in 2^{<\k}$ and $C\subset N_p$ is any co-meager subset of $N_p$, then there
  is no continuous function $C\to 2^\k$ such that $(\eta,\xi)\in E_{S}\cap C^2\iff f(\eta)=f(\xi)$.
\end{Cor}
\begin{proof}
  Apply the same proof as 
  for the first part of Theorem \ref{thm:NotToID}; take $C$ instead of
  $C(f)$ and work inside $N_p$, e.g. instead of $\eta,\bar\eta$ take $p\cat\eta,p\cat\bar\eta$
  for suitable $\eta\in 2^{[\dom p,\k)}$.
\end{proof}

\begin{Def}\label{def:nonrefl} 
  A set $A\subset \k$ \emph{does not reflect} to an ordinal $\a$, if the set $\a\cap A$
  is non-stationary in $\a$, i.e. there exists a closed unbounded subset of $\a$ outside of $A\cap\a$.
\end{Def}

\begin{Thm}\label{thm:nonrefl}
  If $\k=\l^+>\o_1$ and $\square^\k_\mu$ holds, $\mu\le\l$, then for every stationary $S\subset S^\k_\o$,
  there exists a set $B^\mu_{\nr}(S)\subset S$ 
  ($\nr$ for non-reflecting) such that $B^\mu_{\nr}(S)$ does not reflect to any 
  $\a\in S^\k_{\le\mu}\cap S^\k_{\ge\o_1}$ and the sets 
  $\lim C_\a$ witness that, where $\la C_\a\mid \a\in S^\k_{\le\mu}\ra$
  is the $\square_\l$-sequence, i.e. $\lim C_\a\subset \a\setminus B_{\nr}^\mu(S)$ 
  for $\a\in S^\k_{\le\mu}\cap S^\k_{\ge\o_1} $. Since $\cf(\a)>\o$, $\lim C_\a$
  is cub in~$\a$.
\end{Thm}
\begin{proof}
  This is a well known argument and can be found in~\cite{4Jech}. Let $g\colon S\to\k$
  be the function defined by $g(\a)=\OTP(C_\a)$. By the definition of $\square_\mu$, 
  $\OTP(C_\a)<\mu$ for $\a\in S^\k_\o$, so for $\a>\mu$ we have $g(\a)<\a$. By Fodor's 
  lemma there exists a stationary $B_{\nr}^\mu(S)\subset S$ such that
  $\OTP(C_\a)=\OTP(C_\b)$ for all $\a,\b\in B^\mu_{\nr}(\mu)$. If $\a\in \lim C_\b$, then 
  $C_\a=C_\b\cap \a$ and therefore $\OTP(C_\a)<\OTP(C_\b)$. 
  Hence $\lim C_\b\subset \b\setminus B_{\nr}^\mu(S)$.
\end{proof}

\begin{Def}\label{def:Times}
  Let $E_i$ be equivalence relations on $2^{\k\times\{i\}}$ for all $i<\a$ where $\a<\k$.
  Let $E=\bigotimes_{i<\a}E_i$ be an equivalence relation on the space
  $2^{\k\times\a}$ such that $(\eta,\xi)\in E$ if and only if
  for all $i<\a$, $(\eta\rest (\k\times\{i\}),\xi\rest (\k\times\{i\}))\in E_i$.

  Naturally, if $\a=2$, we denote $\bigotimes_{i<2}E_i$ by just $E_0\otimes E_1$
  and we constantly identify $2^{\k\times\{i\}}$ with $2^\k$.
\end{Def}

\begin{Def}\label{def:Plus}
  Given equivalence relations $E_i$ on $2^{\k\times\{i\}}$ for $i<\a<\k^+$, let\label{page:oplus}
  $\bigoplus_{i\in I}E_i$
  be an equivalence relation on $\Cup_{i<\a}2^{\k\times\{i\}}$
  such that $\eta$ and $\xi$ are equivalent if and only if 
  for some $i<\a$, $\eta,\xi\in 2^{\k\times\{i\}}$ and $(\eta,\xi)\in E_i$.

  Intuitively the operation $\oplus$ is taking disjoint unions of the equivalence relations.
  As above, if say $\a=2$, we denote $\bigoplus_{i<2}E_i$ by just $E_0\otimes E_1$
  and we identify $2^{\k\times\{i\}}$ with $2^\k$. 
\end{Def}

\begin{Thm}\label{thm:nonRed} 
  Assume that $\l\in\reg\k$ and $\GC_\l$-characterization holds for~$\k$. 
  \begin{myEnumerate}
  \item \label{nr1} Suppose that 
    $S_1,S_2\subset S^\k_{\ge\l}$ and that $(S_2\setminus S_1)\cap S^\k_\l$ is stationary.
    Then the following holds:
    \begin{myAlphanumerate}
    \item \label{nr1a} $E_{S_1}\not\le_B E_{S_2}$.
    \item \label{nr1b} If $p\in 2^{<\k}$ and $C\subset N_p$ is any co-meager subset of $N_p$, then there
      is no continuous function $C\to 2^\k$ such that $(\eta,\xi)\in E_{S_1}\cap C^2\iff (f(\eta),f(\xi))\in E_{S_2}$.    
    \end{myAlphanumerate}
  \item \label{nr2} Assume that $\k=\l^+>\o_1$, $\mu\in\reg(\k)\setminus\{\o\}$ 
    and $\square^\k_\mu$ holds.
    Let $S\subset S^\k_\o$ be any stationary set and 
    let $B^\mu_\nr(S)$ be the set defined by Theorem \ref{thm:nonrefl}. 
    Then the following holds:
    \begin{myAlphanumerate}
    \item \label{nr2a} Suppose that
      $S_1,S_2\subset S^\k_{\mu}$, $B\subset B_\nr^\mu(S)$ 
      and let $S_1'=S_1\cup B$, $S_2'=S_2\cup B$.
      If $(S_2'\setminus S_1')\cap S^\k_\mu$ is stationary,
      then $E_{S_1'}\not\le_B E_{S_2'}$.
    \item \label{nr2b} Let $S_1$, $S_2$, $B$, $S_1'$ and $S_2'$ be as above.     
      If $(S_2'\setminus S_1')\cap S^\k_\mu$ is stationary,  $p\in 2^{<\k}$ and 
      $C\subset N_p$ is any co-meager subset of $N_p$, then there
      is no continuous function $C\to 2^\k$ such that
      $(\eta,\xi)\in E_{S_1'}\cap C^2\iff (f(\eta),f(\xi))\in E_{S_2'}$.    
    \end{myAlphanumerate}
  \item \label{nr3} 
    Let $S_1,S_2,A_1,A_2\subset S^\k_{\o}$ 
    be either such that $S_2\setminus S_1$ and $A_2\setminus S_1$ are stationary
    \emph{or} such that $S_2\setminus A_1$ and $A_2\setminus A_1$ are stationary.
    Then the following holds:
    \begin{myAlphanumerate}
    \item \label{nr3a} 
      $E_{S_1}\otimes E_{A_1}\not\le_B E_{S_2}\otimes E_{A_2}$.
    \item \label{nr3b} If $C\subset (2^{\k})^2$ (we identify $2^{\k\times 2}$ with $(2^\k)^2$)
      is a set which is 
      co-meager in some $N_r=\{\eta\in (2^\k)^2\mid \eta\rest \dom r=r\}$, $r\in (2^{\a})^2$, $\a<\k$,
      then there is no continuous function $f$ from $C\cap N_r$ to $(2^\k)^2$ such that
      $(\eta,\xi)\in (E_{S_1}\otimes E_{A_1})\cap C^2\iff (f(\eta),f(\xi))\in E_{S_2}\otimes E_{A_2}$.
    \end{myAlphanumerate}
    \item\label{nr4} Assume that $S_1,S_2,A_2\subset \k$ are such that $A_2\setminus S_1$ and
      $S_2\setminus S_1$ are $\o$-stationary. Then
      \begin{myAlphanumerate}
      \item\label{nr4a} $E_{S_1}\not\le_B E_{S_2}\otimes E_{A_2}$.
      \item\label{nr4b} If $p\in 2^{<\k}$ and 
        $C\subset N_p$ is any co-meager subset of $N_p$, then there
        is no continuous function $C\to (2^\k)^2$ such that
        $(\eta,\xi)\in E_{S_1}\cap C^2\iff (f(\eta),f(\xi))\in E_{S_2}\otimes E_{A_2}$.    
      \end{myAlphanumerate}
    \item\label{nr5} Assume that $S_1,A_1,S_2,A_2\subset \k$ are such that 
      $A_2\setminus A_1$ is $\o$-stationary. Then
      \begin{myAlphanumerate}
      \item\label{nr5a} $E_{S_1}\otimes E_{A_1}\not\le_B E_{S_2\cup A_2}$.
      \item\label{nr5b} If $p\in (2^{<\k})^2$ and 
        $C\subset N_p$ is any co-meager subset of $N_p$, then there
        is no continuous function $C\to 2^\k$ such that
        $(\eta,\xi)\in (E_{S_1}\otimes E_{A_1})\cap C^2\iff (f(\eta),f(\xi))\in E_{S_2\cup A_2}$.    
      \end{myAlphanumerate}
  \end{myEnumerate}
\end{Thm}
\begin{proof}
  Item \ref{nr1b} of the theorem implies item \ref{nr1a} as 
  well as all (b)-parts imply the corresponding (a)-parts, because if $f\colon 2^\k\to 2^\k$ 
  is a Borel function, then it is continuous on the co-meager set $C(f)$ 
  (see Definition \ref{def:CoMeagerStuff}).
  Let us start by proving \ref{nr1b}:

  Assume that $S_2\setminus S_1$ is $\l$-stationary, $p\in 2^{<\k}$, $C\subset N_p$ and assume 
  that $f\colon C\to 2^\k$ is a continuous function as described in the Theorem. 
  Let us derive a contradiction.
  Define a strategy for player $\PlTwo$ in the game $\GC_\l(\k\setminus(S_2\setminus S_1))$ as follows.

  Denote the $i$:th move of player $\PlOne$ by $\a_i$ and the $i$:th move of player
  $\PlTwo$ by $\b_i$.
  During the game, at the $i$:th move, $i<\l$, player $\PlTwo$
  secretly defines functions $p_i^0,p_i^1,q_i^0,q_i^1\in 2^{<\k}$
  in such a way that for all $i$ and all $j <  i$ we have
  \begin{myAlphanumerate}
  \item $\dom p_j^0=\dom p_j^1=\b_{j}$ and $\a_{j}\le \dom q_{j+1}^0=\dom q_{j+1}^1\le \a_{j}$,
    and if $j$ is a limit, then $\sup_{i<j}\a_{i}\le \dom q_{j}^0=\dom q_{j}^1\le \b_{j}$,
  \item $p_j^0\subset p_{j+1}^0$, $p_i^1\subset p_{i+1}^1$, $q_i^0\subset q_{i+1}^0$ and $q_i^1\subset q_{i+1}^1$,
  \item $f[C\cap N_{p_i^0}]\subset N_{q_i^0}$ and $f[C\cap N_{p_i^1}]\subset N_{q_i^1}$.
  \end{myAlphanumerate}
  Suppose it is $i$:th move and $i=\g+2k$ for some $k<\o$ and $\g$ which is either $0$ or a limit ordinal,
  and suppose that the players have picked the sequences $(\a_{j})_{j\le i}$ 
  and $(\b_{j})_{j<i}$. Additionally $\PlTwo$ has secretly picked the sequences
  $$(p_i^0)_{i<j},(p_i^1)_{i<j},(q_i^0)_{i<j},(q_i^1)_{i<j}$$
  which satisfy conditions (a)--(c).
  Assume first that $i$ is a successor. 
  If $q^0_{i-1}$ is not $E^*_{S_2}(\dom q^0_{i-1})$-equivalent to $q^1_{i-1}$, then player $\PlTwo$
  plays arbitrarily. Otherwise,
  to decide her next move, player $\PlTwo$ 
  uses Lemma \ref{lem:YouCanFind} (page~\pageref{lem:YouCanFind})
  to find $\eta\in 2^{[\b_{i-1},\k)}$ and $\xi=1-\eta$, 
  such that $p_{i-1}^0\cat\eta\in C$
  and $p_{i-1}^1\cat\xi\in C$. Then she finds $\b_{i}'>\a_i$ such that $f(p_{i-1}^0\cat\eta)(\d)\ne f(p_{i-1}^1\cat\xi)(\d)$
  for some $\d\in [\a_{i},\b_{i}')$. 
  This is possible since $f$ is a reduction and $(q^0_{i-1},q^1_{i-1})\in E^*_{S_2}$.
  Then she picks $\b_{i}>\b_{i}'$ so that
  $$f[C\cap N_{(p_{i-1}^0\cat\eta)\restl\b_{i}}]\subset N_{f(p_{i-1}^0\cat \eta)\restl\b_{i}'}$$
  and 
  $$f[C\cap N_{(p_{i-1}^1\cat\xi)\restl\b_{i}}]\subset N_{f(p_{i-1}^1\cat \xi)\restl\b_{i}'}.$$
  This choice is possible by the continuity of $f$.
  Then she (secretly) sets $p_{i}^0=(p_{i-1}^0\cat\eta)\rest \b_{i}$, $p_i^1=(p_{i-1}^1\cat \xi)\rest\b_{i}$, 
  $q_{i}^0=f(p_{i-1}^0\cat\eta)\rest\b_{i}'$ and $q_{i}^1=f(p_{i-1}^1\cat\xi)\rest \b_{i}'$.
  Note that the new partial functions secretly picked by $\PlTwo$ satisfy conditions (a)--(c).
  
  If $i$ is a limit, then player $\PlTwo$ proceeds as above but instead of $p_{i-1}^n$ she
  uses $\Cup_{i'<i}p_{i'}^n$, $n\in\{0,1\}$, and instead of $\b_{i-1}$ she uses $\sup_{i'<i}\b_{i'}$.
  If $i$ is $0$, then proceed in the same way assuming $p^0_{-1}=p^1_{-1}=q^0_{-1}=q^1_{-1}=\es$
  and $\a_{-1}=\b_{-1}=0$.

  Suppose $i=\g+2k+1$ where $\g$ is again a limit or zero and $k<\o$. 
  Then the moves go in the same way, except that she sets $\eta=\xi$ instead of
  $\eta=1-\xi$ and requires  $f(p_{i-1}^0\cat\eta)(\d)= f(p_{i-1}^1\cat\xi)(\d)$
  for some $\d\in [\a_{i-1},\b_{i}')$ instead of  $f(p_{i-1}^0\cat\eta)(\d)\ne f(p_{i-1}^1\cat\xi)(\d)$
  for some $\d\in [\a_{i-1},\b_{i}')$. Denote this strategy by~$\s$.

  Since $S_2\setminus S_1$ is stationary and $\GC_{\l}$-characterization holds for $\k$, 
  player $\PlOne$ is able play against this strategy such that
  $\sup_{i<\l}\a_i\in S_2\setminus S_1$. Suppose they have played the game to the
  end, so that player $\PlTwo$ used $\sigma$, player $\PlOne$ has won and
  they have picked the sequence $\la\a_i,\b_i\mid i<\l\ra$.
  Let
  $$\a_\l=\sup_{i<\l}\a_i=\sup_{i<\l}\b_i=\sup_{i<\l}\dom p_i=\sup_{i<\l}\dom q_i$$
  and 
  $$p_{\l}^0=\Cup_{i<\l}p_i^0,\ p_{\l}^1=\Cup_{i<\l}p_i^1,\ q_{\l}^0=\Cup_{i<\l}q_i^0\text{ and }q_{\l}^1=\Cup_{i<\l}q_i^1.$$ 
  By continuity, $p_\l^0$, $p_\l^1$, $q_\l^0$ and $q_\l^1$ satisfy condition (c) above
  and $\dom p_\l^0=\dom p_\l^1=\dom q_\l^0=\dom q_\l^1=\sup_{i<\l}\a_i=\sup_{i<\l}\b_i$, so $\a_\l$ is well defined.
  
  On one hand $q_\l^0$ and $q_\l^1$ cannot be extended in an $E_{S_2}$-equivalent way, since
  either they cofinally get same and different values below $\a_\l\in S_2$, or they
  are not $E^*_{S_2}(\g)$-equivalent already for some $\g<\a_\l$. On the other hand
  $p_\l^0$ and $p_\l^1$ can be extended in an $E_{S_1}$-equivalent way, since $\a_\l$ is not in $S_1$ and
  for all $\g<\l$, $\sup_{i<\g}\a_\g$ is not $\mu$-cofinal for any $\mu\ge\l$,
  so cannot be in $S_1$ either~$(*)$. 

  Let $\eta,\xi\in 2^\k$ be extensions of $p_\l^0$ and $p_{\l}^1$ respectively such that $(\eta,\xi)\in E_{S_1}\cap C^2$.
  Now $f(\eta)$ and $f(\xi)$ cannot be $E_{S_2}$-equivalent, since by condition (c),
  they must extend $q^0_\l$ and $q^1_\l$ respectively.

  Now let us prove \ref{nr2b} which implies \ref{nr2a}. 
  Let \mbox{$\la C_\a^\mu\mid \a\in S^\k_{\le\mu}\ra$} 
  be the $\square^\k_\mu$-sequence and denote by $t^\mu$ the function
  $\a\mapsto C_\a^\mu$.

  Let player $\PlTwo$ define her strategy in the game $\GC(\k\setminus (S_2'\setminus S_1'))$
  exactly as in the proof of \ref{nr1b}. Note that $S_2'\setminus S_1'=S_2\setminus S_1$
  since $\mu>\o$.
  Denote this strategy by $\sigma$. We know that, as above, Player $\PlOne$
  is able to beat $\sigma$. However, now it is not enough, because in order to be able to
  extend $p_\mu^0$ and $p_\mu^{1}$ in an $E_{S_1'}$-equivalent way, he needs to ensure that
  $$S_1'\cap \lim{}\!_{\o}(\{\a_i\mid i<\mu\})=\es\eqno(**)$$
  where $\lim_\o X$ is the set of $\o$-limits of elements of $X$, i.e. we cannot rely on the sentence
  followed by~$(*)$~above. On the other hand $(**)$ is 
  sufficient, because $S'_1\subset S^\k_\mu\cup S^\k_{\o}$.
  
  Let us show that it is possible for player $\PlOne$ to play against $\sigma$ as required.

  Let $\nu>\k$ be a sufficiently large cardinal and let $M$ be an elementary submodel of
  $\la H_{\nu}, \sigma, \k, t^\mu \ra$ such that $|M|<\k$ and $\a=\k\cap M$ is an ordinal in 
  $S_2'\setminus S_1'$. 

  In the game, suppose that the sequence $d=\la \a_j,\b_j\mid j<i\ra$ has been played
  before move $i$ and suppose that this sequence is in $M$. Player $\PlOne$ will
  now pick $\a_i$ to be the smallest element in $C_\a^\mu$ which is above
  $\sup_{j<i}\b_j$. Since $C_\a^\mu\cap \b=C_\b^\mu$ for any $\b\in \lim C_\a^\mu$ and $C_\b^\mu\in M$,
  this element is definable in $M$ from the sequence~$d$ and $t^\mu$. This guarantees
  that the sequence obtained on the following move is also in $M$. At limits
  the sequence is in $M$, because it is definable from $t^\mu$ and $\sigma$. 
  Since $\OTP(C_\a^\mu)=\mu$, the game ends at 
  $\a$ and player $\PlOne$ wins. Also the requirement $(**)$ is satisfied because
  he picked elements only from $C_\a^\mu$ and 
  so $\lim_\o\{\a_i\mid i<\mu\}\subset\lim_\o(C^\mu_\a)\subset \a\setminus B$
  which gives the result.

  Next let us prove \ref{nr3b} which again implies \ref{nr3a}. 
  The proofs of \ref{nr4} and \ref{nr5} are very similar to that of \ref{nr3} 
  and are left to the reader.

  So, let $S_1$, $A_1$, $S_2$, $A_2$, $C$ 
  and $r$ be as in the statement of \ref{nr3} and suppose that
  there is a counter example~$f$. Assume that $S_2\setminus S_1$ and $A_2\setminus S_1$ are stationary,
  the other case being symmetric.
  Let us define the property $P$:
  \begin{myItemize}
  \item[$P$:] There exist $p,p'\in (2^\a)^2$, $p=(p_1,p_2)$ and $p'=(p'_1,p'_2)$,
    such that 
    \begin{myAlphanumerate}
    \item $r\subset p\cap p'$, 
    \item $p_2=p'_2$, $(p_1,p'_1)\in E^*_{S_1}(\a+1)$ (see Definition~\ref{def:EqRels}, page~\pageref{def:EqRels}),      
    \item for all $\eta\in C\cap N_p$ and $\eta'\in C\cap N_{p'}$, $\eta=(\eta_1,\eta_2)$,
      $\eta'=(\eta'_1,\eta'_2)$, if
      $\eta_2=\eta'_2$ and $(\eta_1,\eta'_1)\in E^*_{S_1}$, then
      $f(\eta)_1 \sd f(\eta')_1\subset \dom p_1$ where
      $f(\eta)=(f(\eta)_1,f(\eta)_2)$.
    \end{myAlphanumerate} 
  \end{myItemize}
  We will show that both $P$ and $\lnot P$ lead to a contradiction.
  Assume first $\lnot P$. Now the argument is similar to the proof
  of \ref{nr1b}. Player $\PlTwo$ defines her strategy in the same way
  but this time she chooses the elements $p^n_i$ and $q^n_i$ from $(2^{\a})^2$
  instead of $2^\a$ so that $p^n_i=(p^n_{i,1},p^n_{i,2})$, $q^n_i=(q^n_{i,1},q^n_{i,2})$
  and for all $i<\l$, $p^0_{i,2}=p^1_{i,2}$. 
  In building the strategy she looks only at $q^n_{i,1}$ and ignores $q^n_{i,2}$. In
  other words she pretends that the game is for $E_{S_1}$ and $E_{S_2}$ in the proof of \ref{nr1}.
  At the even moves she extends $p^0_{i,1}$ and $p^1_{i,1}$ by $\eta$ and $\eta'$
  which witness the failure of item (c) (but not of (a) and (b)) 
  of property $P$ for $p^0_{i}$ and $p^1_{i}$.
  Then there is $\a\in f(\eta)_1\sd f(\eta')_1$, $\a>\dom p^0_{i,1}$. And then she chooses 
  $q^0_{i,1}$ and $q^1_{i,1}$ to be initial segments of $f(\eta)_1$ and $f(\eta')_1$ respectively.
  
  At the odd moves she just extends $p^0_{i,1}$ and $p^1_{i,1}$ in an $E_{S_1}$-equivalent way,
  so that she finds an $\a>\dom p^0_{i,1}$, $q^0_{i,1}$ and $q^1_{i,1}$ such that $q^0_{i,1}(\a)=q^1_{i,1}(\a)$
  and $f[N_{p^0_i}\cap C]\subset N_{q^0_i}$. 
  
  As in the proof of \ref{nr1}, $\PlOne$ responses by playing towards an ordinal in $S_2\setminus S_1$.
  During the game they either hit a point at which $q^0_{i,2}$ and $q^1_{i,2}$ cannot be extended
  to be $E_{A_2}$-equivalent or else they play the game to the end whence 
  $q^0_{\l,1}$ and $q^1_{\l,1}$ cannot be extended in a $E_{S_2}$-equivalent way
  but $p^0_\l$ and $p^1_\l$ can be extended to $E_{S_1}\otimes E_{A_1}$-equivalent way.

  Assume that $P$ holds. Fix $p$ and $p'$ which witness that.
  Now player $\PlTwo$ builds her strategy as if they were playing 
  between $E_{S_1}$ and $E_{A_2}$. This time she concentrates on $q^0_{i,2}$ and $q^1_{i,2}$ 
  instead of $q^0_{i,1}$ and $q^1_{i,1}$.
  At the even moves she extends $p^0_{i,1}$ and $p^1_{i,1}$ by $\eta$ and $\bar\eta$ respectively
  for some $\eta$. Also, as above, $p^0_{i,2}$ and $p^1_{i,2}$
  are extended in the same way. By item (c) $f(\eta)_1\sd f(\eta')_1$ is bounded by $\dom p^0_{i,1}$,
  but $f(\eta)$ and $f(\eta')$ can't be $E_{S_2}\otimes E_{A_2}$-equivalent, because $f$ is assumed to be a reduction.
  Hence there must exist $\a>\dom p^0_{i,1}$, $q^0_{i,2}$ and $q^0_{i,2}$ such that $q^0_{i,2}(\a)\ne q^1_{i,2}(\a)$.
  The rest of the argument goes similarly as above.
\end{proof}

\begin{Cor}\label{cor:E_0nonredtoE_S}
  If $\GC_{\l}$-characterization holds for $\k$ and $S\subset \k$ is $\l$-stationary, 
  then $E_0\not\le E_S$.
  In particular, if $S$ is $\o$-stationary, then $E_0\not\le E_S$.
\end{Cor}
\begin{proof}
  Follows from Theorem \ref{thm:nonRed}.\ref{nr1a} by taking $S_1=\es$, since 
  $E_\es=E_0$ and $\GC_{\o}$-characterization holds for $\k$.
\end{proof}

\begin{Cor}\label{cor:Antichain}
  There is an antichain\footnote{By an antichain I refer here to a 
    family of pairwise incomparable elements unlike e.g. in forcing context.} 
  of Borel equivalence relations on $2^\k$ of length~$2^\k$.
\end{Cor}
\begin{proof}
  Take disjoint $\o$-stationary sets $S_i$, $i<\k$. 
  Let $f\colon \k\times 2\to\k$ be a bijection. For each $\eta\in 2^\k$ let
  $A_\eta=\{(\a,n)\in \k\times 2\mid (n=0\land \eta(\a)=1)\lor(n=1\land \eta(\a)=0)\}$. 
  For each $\eta\ne\xi$ clearly
  $A_\eta\setminus A_\xi\ne\es\ne A_\xi\setminus A_\eta$. Let 
  $$S_\eta=\Cup_{i\in f[A_\eta]} S_i.$$
  Now $\{E_{S_\eta}\mid \eta\in 2^\k\}$ is an antichain by Theorem~\ref{thm:nonRed}.\ref{nr1b}.
\end{proof}


Let us show that all these relations are below $E_0$. It is already shown that they are not above it 
(Corollary~\ref{cor:E_0nonredtoE_S}), provided $\GC_\l$-characterization holds for~$\k$.
Again, similar ideas will be used in the proof of Theorems~\ref{thm:Main1} and~\ref{thm:Main2}.

\begin{Thm}\label{thm:E_SEmbedstoE_0}
  For all $S$, $E_S\le_B E_0$.
\end{Thm}
\begin{proof}
  Let us show that $E_S$ is reducible to $E_0$ on $\k^\k$ which is in turn bireducible with $E_0$ on $2^\k$
  (see \cite{5FriHytKul}).
  Let us define an equivalence relation $\sim$ on $2^{<\k}$ as on page
  \pageref{pagesimrel}, such that $p\sim q$ if and only if 
  $\dom p=\dom q$ and $p\sd q$ is eventually
  constant, i.e. for some $\a<\dom p$, $(p\sd q)(\g)$ is the same for all $\g\in [\a,\dom p)$.
  Let $s\colon 2^{<\k}\to \k$ be a map such that $p\sim q\iff s(p)=s(q)$. 
  Let $\{A_i\mid i\in S\}$ be a partition of $\lim \k$ into disjoint unbounded sets.
  Suppose $\eta\in 2^\k$ and
  define $f(\eta)=\xi\in \k^\k$ as follows.
  \begin{myItemize}
  \item If $\a$ is a successor, $\a=\b+1$, then $\xi(\a)=\eta(\b)$.
  \item If $\a$ is a limit, then $\a\in A_i$ for some $i\in S$. Let $\xi(\a)=s(\eta\rest i)$
  \end{myItemize}
  Let us show that $f$ is the desired reduction from $E_S$ to $E_0$. 
  Assume that $\eta$ and $\xi$ are
  $E_S$-equivalent. If $\a$ is a limit and $\a\in A_i$, then, since $\eta$ and $\xi$ are $E_S$-equivalent,
  we have $\eta\rest i\sim\xi\rest i$, so $s(\eta\rest i)=s(\xi\rest i)$ and so $f(\eta)(\a)=f(\xi)(\a)$.
  There is $\b$ such that $\eta(\g)=\xi(\g)$ for all $\g>\b$. This implies that
  for all successors $\g>\b$ we also have $f(\eta)(\g)=f(\xi)(\g)$. Hence $f(\eta)$ and $f(\xi)$ are $E_0$-equivalent.
  Assume now that $\eta$ and $\xi$ are not $E_S$-equivalent. Then there are two cases: 
  \begin{enumerate}
  \item $\eta\sd\xi$ is unbounded. Now $f(\eta)(\b+1)=\eta(\b)$ and $f(\xi)(\b+1)=\xi(\b)$ for all $\b$, so we have
    $$\{\b\mid \eta(\b)\ne\xi(\b)\} =\{\b\mid f(\eta)(\b+1)\ne\xi(\b+1)\}.$$
    If the former is unbounded, then so is the latter.
  \item For some $i\in S$, $\eta\rest i\not\sim \xi\rest i$. This implies that 
    $f(\eta)(\a)\ne f(\xi)(\a)$ for all $\a\in A_i$.
    and we get that $\{\b\mid f(\eta)(\b)\ne\xi(\b)\}$ is again unbounded. 
  \end{enumerate}
  It is easy to check that $f$ is continuous.
\end{proof}



\subsection{Proofs of the Main Theorems}

\begin{proof}[Proof of Theorem \ref{thm:Main1}]\label{pageproof1}
  The subject of the proof is that for a regular $\l<\k$, if $\GC_\l$-characterization holds for $\k$, then 
  the order $\la\Po(\k),\subset_{\NS(\l)}\ra$ can be embedded into $\la \E^B_\k,\le_B\ra$ strictly below $E_0$
  and above $\id_{2^\k}$. 


  Let $h\colon \o\times\k\to\k$ be a bijection. 
  Let $\tilde h\colon 2^{\o\times\k}\to 2^\k$ be defined by $\tilde h(\eta)(\a)=\eta(h^{-1}(\a))$.
  We define the topology on $2^{\o\times\k}$
  to be generated by the sets $\{\tilde h^{-1}V\mid V\text{ is open in }2^\k\}$. Then $\tilde h$ is a homeomorphism
  between $2^{\o\times\k}$ and $2^\k$.
  If $g\colon\k\times\k\to\k$ is a bijection, we similarly get a topology onto $2^{\k\times\k}$ and
  a homeomorphism $\tilde g$ from $2^{\k\times\k}$ onto $2^\k$. By combining these two we get a homeomorphism between
  $2^{\o\times\k}\times 2^\k$ and $2^\k$, and so without loss of generality we 
  can consider equivalence relations on these spaces.

  For a given equivalence relation $E$ on $2^\k$, let $\overline{E}$ be the equivalence relation on
  $2^{\o\times\k}\times 2^\k$ defined by 
  $$((\eta,\xi),(\eta',\xi'))\in \overline{E}\iff \eta=\eta'\land (\xi,\xi')\in E.$$
  Essentially $\overline{E}$ is the same as $\id\otimes E$, since $2^{\o\times\k}\approx 2^\k$.

  \begin{RemarkN}\label{rem:nonRedbar}
    Corollary \ref{cor:NotToID2}, Theorem \ref{thm:nonRed} and Corollary \ref{cor:E_0nonredtoE_S} hold
    even if $E_{S}$ is replaced everywhere by $\overline{E_S}$ for all $S\subset\k$.
    \begin{proof}
      Let us show this for Theorem \ref{thm:nonRed}.\ref{nr1}.
      The proof goes exactly as the proof of Theorem \ref{thm:nonRed}.\ref{nr1}, but
      player $\PlOne$ now picks the functions $p^n_k$ from $\Cup_{\a<\k}2^{\o\times\a}\times 2^\a$
      instead of $2^{<\k}$, $p^n_k=(p^n_{k,1},p^n_{k,2})$, and 
      requires that at each move $p^0_{k,1}=p^1_{k,1}$. Otherwise the argument proceeds in the
      same manner. Similarly for \ref{thm:nonRed}.\ref{nr2}, \ref{thm:nonRed}.\ref{nr3}, 
      \ref{thm:nonRed}.\ref{nr4} and \ref{thm:nonRed}.\ref{nr5}.

      Modify the proof of the first part of Theorem \ref{thm:NotToID} in a similar way 
      to obtain the result for
      Corollary~\ref{cor:NotToID2}. Corollary~\ref{cor:E_0nonredtoE_S} follows from the 
      modified version of Theorem~\ref{thm:nonRed}.
    \end{proof}
  \end{RemarkN}
  For $S\subset \k$ let 
  $$G(S)=\overline{E_{S^\k_\l\setminus S}}.$$

  Let us show that $G\colon \Po(\k)\to\E^{B}_{\k}$ is the desired embedding.
  Without loss of generality let us assume that $G$ is restricted to $\Po(S^\k_\l)$, whence stationary is the same
  as $\l$-stationary and non-stationary is the same as not $\l$-stationary.
  For arbitrary $S_1,S_2\subset S^\k_\l$ we have to show:
  \begin{enumerate}
  \item If $S_2\setminus S_1$ is stationary, then $\overline{E_{S_1}}\not\le_B \overline{E_{S_2}}$
  \item If $S_2\setminus S_1$ is non-stationary, then $\overline{E_{S_1}}\le_B \overline{E_{S_2}}$
  \item $\id_{2^\k}\lneq_B\overline{E_{S_1}}\lneq_B E_0$.
  \end{enumerate}
  If $\eta\in 2^{\o\times\k}$, denote $\eta_i(\a)=\eta(i,\a)$ and $( \eta_i)_{i<\o}=\eta$.
  \begin{claim}{1}
    If $S_2\setminus S_1$ is stationary, then $\overline{E_{S_1}}\not\le_B \overline{E_{S_2}}$. 
    Also $E_0\not\le \overline{E_S}$.
  \end{claim}
  \begin{proof}
    Follows from Theorem \ref{thm:nonRed}.\ref{nr1a} and Remark \ref{rem:nonRedbar}.
  \end{proof}
  \begin{claim}{2}\label{pageclaimtwo}
    If $S_2\setminus S_1$ is non-stationary,
    then $\overline{E_{S_1}}\le_B \overline{E_{S_2}}$.
  \end{claim}
  \begin{proof}
    Let us split this into two parts accorind to the stationarity of $S_2$.
    Assume first that $S_2$ is non-stationary. Let $C$ be a cub set outside $S_2$.
    Let $f\colon 2^\k\to 2^{\o\times\k}\times 2^\k$ be the function defined as follows.
    For $\eta\in 2^\k$ let $f(\eta)=\la (\eta_i)_{i<\o},\xi\ra$ be such that $\eta_i(\a)=0$
    for all $\a<\k$ and $i<\o$ and $\xi(\a)=0$ for all $\a\notin C$. If $\a\in C$, then 
    let $\xi(\a)=\eta(\OTP(\a\cap C))$. This is easily verified to be a reduction 
    from $E_0$ to $\overline{E_{S_2}}$. By the following Claim 3, $\overline{E_{S_1}}\le_B E_0$,
    so we are done. 

    Assume now that $S_2$ is stationary. Note that then $S_1$ is also stationary.
    Let $C$ be a cub set such that $S_2\cap C\subset S_1$.
    Assume that $\la (\eta_i)_{i<\o},\xi\ra\in 2^{\o\times \k}\times 2^\k$ and let us define
    $$f(\la (\eta_i)_{i<\o},\xi\ra)=\la (\eta'_i)_{i<\o},\xi'\ra\in 2^{\o\times \k}\times 2^\k$$
    as follows. 
    For $i\ge 0$ let
    $$\eta'_{i+1}=\eta_i.$$
    For all $\a<\k$, let $\xi'(\a)=\xi(\min(C\setminus \a))$. Then let $s$ be the function as defined
    in the proof of Theorem \ref{thm:NotToID} (on page \pageref{pagesimrel}) 
    and for all $\a<\k$ let $\b(\a)$ be the $\a$:th element of
    $S_1\setminus S_2$. For all $\a<\k$, let
    $$\eta'_0(\a)=s(\xi\rest\b(\a)).$$
    Let us show that this defines a continuous reduction. 

    Suppose $\la (\eta^0_i)_{i<\o},\xi^0\ra$ and $\la (\eta^1_i)_{i<\o},\xi^1\ra$ are $\overline{E_{S_1}}$-equivalent.
    Denote their images under $f$ by $\la (\rho^0_i)_{i<\o},\zeta^0\ra$ and  $\la (\rho^1_i)_{i<\o},\zeta^1\ra$ 
    respectively. Since $\eta^0_i=\eta^1_i$ for all $i<\o$, we have $\rho^0_i=\rho^1_i$ for all $0<i<\o$.
    Since for all $\a\in S_1$ we have that $\xi^0\rest\a$ and $\xi^1\rest \a$ are $\sim$-equivalent 
    (as in the definition of $s$),
    we have that $\rho_0^0(\b)=\rho_0^1(\b)$ for all $\b<\k$. 

    Suppose now that $\a\in S_2$. The aim is to show that $\zeta^0\rest \a\sim \zeta^1\rest\a$.
    If $\a\notin C$, then there is $\b<\a$ such that
    $C\cap (\b,\a)=\es$, because $C$ is closed. This implies that
    for all $\b<\g<\g'<\a$, $\min(C\setminus \g')=\min(C\setminus\g)$, 
    so by the definition of $f$, $\zeta^0(\g)=\zeta^0(\g')$ and $\zeta^1(\g)=\zeta^1(\g')$. 
    Now by fixing $\g_0$ between $\b$ and $\a$ we deduce that $\zeta^0\rest(\b,\a)$ is constant and $\zeta^1\rest(\b,\a)$
    is constant, since for all $\g<\a$ we have $\zeta^0(\g)=\zeta^0(\g_0)$ and $\zeta^1(\g)=\zeta^1(\g_0)=\zeta^1(\g)$. 
    Hence $(\zeta^0\sd\zeta^1)\rest (\b,\a)$ is constant which by the definition of $\sim$ implies that $\zeta^0\rest\a\sim \zeta^1\rest\a$.

    If $\a\in C$, then, since $\a$ is also in $S_2$, we have by the definition of $C$ that $\a\in S_1$.
    Thus, there is $\b<\a$ such that $(\xi^0\sd \xi^1)\rest (\b,\a)$ is constant which implies 
    that for some $k\in \{0,1\}$ we have $(\zeta^0\sd \zeta^1)(\g)=k$ for all $\g\in (\b,\a)\cap C$. But if $\g\in (\b,\a)\setminus C$,
    then, again by the definition of $f$, we have $(\zeta^0\sd \zeta^1)(\g)=(\zeta^0\sd \zeta^1)(\g')$ for some $\g\in (\b,\a)\cap C$,
    so $(\zeta^0\sd \zeta^1)(\g)$ also equals to $k$.

    This shows that $\zeta^0$ and $\zeta^1$ are $E_{S_2}^*$-equivalent. 
    It remains to show that they are $E_0$-equivalent. But
    since $\xi^0$ and $\xi^1$ are $E_0$-equivalent, the number $k\in \{0,1\}$ referred above equals $0$ for all $\a$ large enough and
    we are done.

    Next let us show that if $\la (\eta^0_i)_{i<\o},\xi^0\ra$ and $\la (\eta^1_i)_{i<\o},\xi^1\ra$ are 
    not $\overline{E_{S_1}}$-equivalent,
    then $\la (\rho^0_i)_{i<\o},\zeta^0\ra$ and $\la (\rho^1_i)_{i<\o},\zeta^1\ra$ are not $\overline{E_{S_2}}$-equivalent.
    This is just reversing implications of the above argument. If $\eta_i^0\ne\eta_i^1$ for some $i<\o$, then
    $\rho_{i+1}^0\ne\rho_{i+1}^1$, so we can assume that $(\xi^0,\xi^1)\notin E_{S_1}$. If $\xi^0$ and $\xi^1$ are not $E_{S_1}^*$-equivalent,
    then $\rho^0(\a)\ne\rho^1(\a)$ for some $\a<\k$. 

    The remaining case is that $\xi^0$ and $\xi^1$ are $E_{S_1}^*$-equivalent but not
    $E_0$-equivalent.
    But then in fact $\xi^0\sd \xi^1$ is eventually equal to $1$, since otherwise the sets
    $$C_1=\{\a\mid \{\b<\a\mid (\xi^0\sd \xi^1)(\b)=1\}\text{ is unbounded in }\a\}$$
    and
    $$C_2=\{\a\mid \{\b<\a\mid (\xi^0\sd \xi^1)(\b)=0\}\text{ is unbounded in }\a\}$$
    are both cub and by the stationarity of $S_1$,
    there exists a point $\a\in C_1\cap C_2\cap S_1$ which contradicts the fact that $\xi_0$ and $\xi_1$ are
    $E_{S_1}^*$-equivalent. So $\xi^0\sd \xi^1$ is eventually equal to $1$ and this finally implies that also $\zeta^0$ and
    $\zeta^1$ cannot be $E_0$-equivalent.
  \end{proof} 
  \begin{claim}{3}
    Let $S\subset S^\k_\l$. Then $\id\lneq_B \overline{E_S}<_B E_0$.
    If $S$ is stationary, then also $E_0\not\le_B \overline{E_S}$.
  \end{claim}
  \begin{proof}
    If $\eta\in 2^\k$, let $\eta_0=\eta$ and $\eta_i(\a)=\xi(\a)=0$ for all $\a<\k$. Then
    $\eta\mapsto \la(\eta_i)_{i<\o},\xi\ra$ defines a reduction from $\id$ to $\overline{E_S}$. 
    On the other hand $\overline{E_S}$ is not reducible to $\id$ by Remark~\ref{rem:nonRedbar}.

    Let $u\colon 2^{\o\times \k}\to 2^\k$ be a reduction from $\id_{2^{\o\times\k}}$ to $E_0$. 
    Let $v\colon 2^{\k}\to 2^{\k}$ be a reduction from $E_{S}$ to $E_0$ which exists by \ref{thm:E_SEmbedstoE_0}.
    Let $\{A,B\}$ be a partition of $\k$ into two disjoint unbounded subsets. Let $(\eta,\eta')\in 2^{\o\times\k}\times 2^\k$
    and let us define $\xi=f(\eta,\eta')\in 2^\k$. If $\a\in A$, then let $\xi(\a)=u(\eta)(\OTP(\a\cap A))$.
    If $\a\in B$, then let $\xi(\a)=v(\eta')(\OTP(\a\cap B))$. (See page \pageref{page:OTP} for notation.)

    Now if $((\eta_0,\eta'_0),(\eta_1,\eta_1'))\in (2^{\o\times\k}\times 2^\k)^2$ are $\overline{E_S}$-equivalent, then
    both $u(\eta_0)\sd u(\eta_1)$ and $v(\eta_0')\sd v(\eta_1')$ are eventually equal to zero which clearly implies that
    $f(\eta_0,\eta_0')\sd f(\eta_1,\eta_1')$ is eventually zero, and so $f(\eta_0,\eta_0')$ and $f(\eta_1,\eta_1')$ 
    are $E_0$-equivalent. Similarly, if  $(\eta_0,\eta'_0)$ and $(\eta_1,\eta_1')$ are not $\overline{E_S}$-equivalent, then
    either $u(\eta_0)\sd u(\eta_1)$ or $v(\eta_0')\sd v(\eta_1')$ is not eventually zero, and so
    $f(\eta_0,\eta_0')$ and $f(\eta_1,\eta_1')$ are not $E_0$-equivalent.    

    If $S$ is stationary, then $E_0\not\le_B \overline{E_S}$
    by Corollary \ref{cor:E_0nonredtoE_S} and Remark \ref{rem:nonRedbar}.
  \end{proof}\qedhere
\end{proof}

\begin{proof}[Proof of Theorem \ref{thm:Main2}] \label{pageproof2}
  Let us review the statement of the Theorem:
  assuming $\k=\o_1$, or $\k=\l^+$ and $\square_\l$, the
  partial order $\la \Po(\k),\subset_{\NS}\ra$ can be embedded into \mbox{$\la \E^B_\k,\le_B\ra$.}

  If $\k=\o_1$, then this is just the second part (a special case) 
  of Corollary~\ref{cor:Second} on page~\pageref{cor:Second} 
  and follows from Theorem~\ref{thm:Main1}. 


  Recall Definition \ref{def:Plus} on page~\pageref{def:Plus}. Let us see that if $\a<\k$, then
  $\Cup_{i<\a}2^{\k\times\{i\}}$ is homeomorphic to $2^\k$ and so the domains of
  the forthcoming equivalence relations can be thought without loss of generality
  to be $2^\k$. So fix $\a<\k$. For all $\b+1<\a$
  let $\zeta_\b\colon \b+1\to 2$ be the function $\zeta_{\b}(\g)=0$ for all
  $\g<\b$ and $\zeta_\b(\b)=1$ and let 
  $\zeta_\a\colon \a\to 2$ be the constant function with value $0$. 
  Clearly $(\zeta_\b)_{\b\le\a}$ is a maximal antichain. By rearranging the indexation
  we can assume that $(\zeta_\b)_{\b<\a}$ is a maximal antichain.
  If $\eta\in 2^{\k\times\{i\}}$, $i<\a$, let $\xi=\eta+i$ be the function with $\dom\xi=[i+1,\k)$
  and $\xi(\g)=\eta(\OTP(\g\setminus i))$ and let
  $$f(\eta)=\zeta_{i}\cat(\eta+i).$$
  Then $f$ is a homeomorphism $\Cup_{i<\a}2^{\k\times\{i\}}\to 2^\k$.
  

  Assume $S\subset \k$ and let us construct the equivalence relation $H_S$.
  Denote for short $r=\reg\k$, the set of regular cardinals below $\k$. 
  Since $\k$ is not inaccessible, $|r|<\k$. Let $\{K_\mu\subset S^\k_\o\mid \mu\in r\}$ 
  be a partition of $S^\k_\o$ into disjoint stationary sets. For each $\mu\in r\setminus \{\o\}$,
  let $A_\mu=B^\mu_{\nr}(K_\mu)$ be the set given by 
  Theorem~\ref{thm:nonrefl}.
  Additionally let $\{A^0_\o,A^1_\o,A^2_\o,A^3_\o\}$ be a partition of $K_\o$ into disjoint stationary sets.

  Let
  \begin{eqnarray*}\label{page:equationthatrequiersexplanation}
    H_S&=&\phantom{\oplus}(\id_{2^\k}\otimes E_{A^2_\o\cup ((S\cap S^\k_\o)\setminus A^0_\o)}\otimes E_{A^0_\o})\\
    &&\oplus(\id_{2^\k}\otimes E_{A^3_\o\cup ((S\cap S^\k_\o)\setminus A^1_\o)}\otimes E_{A^1_\o})\\
    &&\oplus\bigoplus_{{\mu\in r}\atop{\mu>\o}} (\id_{2^\k}\otimes E_{(S\cap S^\k_\mu)\cup A_\mu}). 
  \end{eqnarray*}
  This might require a bit of explanation. $H_S$ is a disjoint union of the 
  equivalence relations listed in the equation.
  The final part of the equation lists all the relations obtained by
  splitting the set $S$ into pieces of fixed uncountable cofinality and coupling them with 
  the non-reflecting $\o$-stationary sets $A_\mu$.
  The operation $E\mapsto \id_{2^\k}\otimes E$ is the same as the operation $E\mapsto\overline E$ 
  in the proof of Theorem~\ref{thm:Main1} above after the identification
  $2^{\o\times\k}\approx 2^\k$.  The first two lines of the equation
  deal with the $\o$-cofinal part of $S$. It is trickier, because the ``coding sets'' $A_\mu$ also
  consist of $\o$-cofinal ordinals. The way we have built up the relations makes it possible to
  use Theorem~\ref{thm:nonRed} to prove that $S\mapsto H_{\k\setminus S}$ is the desired embedding.

  In order to make the sequel a bit more readable, let us denote
  \begin{eqnarray*}
    \B_\o^0(S)&=&(\id_{2^\k}\otimes E_{A^2_\o\cup ((S\cap S^\k_\o)\setminus A^0_\o)}\otimes E_{A^0_\o}),\\
    \B_\o^1(S)&=&(\id_{2^\k}\otimes E_{A^3_\o\cup ((S\cap S^\k_\o)\setminus A^1_\o)}\otimes E_{A^1_\o}),\\
    \B_\mu(S) &=&(\id_{2^\k}\otimes E_{(S\cap S^\k_\mu)\cup A_\mu}),
  \end{eqnarray*}
  for $\mu\in r\setminus\{\o\}$. With this notation we have
  $$H_S=\B^0_\o(S)\oplus \B^1_\o(S)\oplus\bigoplus_{{\mu\in r}\atop{\mu>\o}} \B_\mu(S).$$


  Let us show that $S\mapsto H_{\k\setminus S}$ is an embedding 
  from $\la\Po(\k),\subset_{\NS}\ra$ into $\la \E^B_\k,\le_B\ra$.
  Suppose $S_2\setminus S_1$ is non-stationary. Then for each $\mu\in r\setminus\{\o\}$ the set 
  $$\big((S^\k_\mu\cap S_2)\cup A_\mu\big)\setminus \big((S^\k_\mu\cap S_1)\cup A_\mu\big)$$
  is non-stationary as well as are the sets
  $$\big(A^2_\o\cup((S_2\cap S^\k_\o)\setminus A^0_\o)\big)\setminus\big(A^2_\o\cup((S_1\cap S^\k_\o)\setminus A^0_\o)\big)$$
  and
  $$\big(A^3_\o\cup((S_2\cap S^\k_\o)\setminus A^1_\o)\big)\setminus\big(A^3_\o\cup((S_1\cap S^\k_\o)\setminus A^1_\o)\big)$$
  so by Claim 2 of the proof of Theorem~\ref{thm:Main1} 
  (page~\pageref{pageclaimtwo}) we have for all $\mu\in r\setminus\{\o\}$ that
  $$(\id_{2^\k}\otimes E_{(S_1\cap S^\k_\mu)\cup A_\mu})\le_B(\id_{2^\k}\otimes E_{(S_2\cap S^\k_\mu)\cup A_\mu}),$$ 
  $$(\id_{2^\k}\otimes E_{A^2_\o\cup ((S_1\cap S^\k_\o)\setminus A^0_\o)})\le_B (\id_{2^\k}\otimes E_{A^2_\o\cup ((S_2\cap S^\k_\o)\setminus A^0_\o)}),$$
  and 
  $$(\id_{2^\k}\otimes E_{A^3_\o\cup ((S_1\cap S^\k_\o)\setminus A^1_\o)})\le_B (\id_{2^\k}\otimes E_{A^3_\o\cup ((S_2\cap S^\k_\o)\setminus A^1_\o)}).$$
  Of course this implies that for all $\mu\in r\setminus\{\o\}$ 
  $$(\id_{2^\k}\otimes E_{A^2_\o\cup ((S_1\cap S^\k_\o)\setminus A^0_\o)}\otimes E_{A^0_\o})\le_B(\id_{2^\k}\otimes E_{A^2_\o\cup ((S_2\cap S^\k_\o)\setminus A^0_\o)}\otimes E_{A^0_\o})$$
  and  that
  $$(\id_{2^\k}\otimes E_{A^3_\o\cup ((S_1\cap S^\k_\o)\setminus A^1_\o)}\otimes E_{A^1_\o})\le_B (\id_{2^\k}\otimes E_{A^3_\o\cup ((S_2\cap S^\k_\o)\setminus A^1_\o)}\otimes E_{A^1_\o})$$
  which precisely means that $\B^0_\o(S_1)\le_B \B^0_\o(S_2)$,
  $\B^1_\o(S_1)\le_B \B^1_\o(S_2)$ and $\B_\mu(S_1)\le_B \B_\mu(S_2)$ for all $\mu\in r\setminus\{\o\}$.
  Combining these reductions we get a reduction from $H_{S_1}$ to $H_{S_2}$.

  Assume that $S_2\setminus S_1$ is stationary. We want to show that $H_{S_1}\not\le_B H_{S_2}$.
  $H_{S_1}$ is a disjoint union the equivalence relations 
  $\B^0_\o(S_1)$, $\B^1_\o(S_1)$ and $\B_\mu(S_1)$ for $\mu\in r\setminus\{\o\}$ .
  Let us call these equivalence relations \emph{the building blocks of $H_{S_1}$} and 
  similarly for $H_{S_2}$.

  Each building block of $H_{S_1}$ can be easily reduced to $H_{S_1}$ via inclusion, 
  so it is sufficient to show that there is one block that
  cannot be reduced to $H_{S_2}$. We will show that if $\mu_1$ is the least cardinal such that
  $S^\k_{\mu_1}\cap(S_2\setminus S_1)$ is stationary, then 
  \begin{myItemize}
  \item that building block is $\B_{\mu_1}(S_1)$, if $\mu_1>\o$.    
  \item that building block is either $\B^0_\o(S_1)$ or $\B^1_\o(S_1)$, if $\mu_1=\o$.
  \end{myItemize}
  Such a cardinal $\mu_1$ exists because $\k$ is not inaccessible and $|r|<\k$.

  Suppose that $f$ is a reduction from a building block of $H_{S_1}$, call it $\B$,
  to $H_{S_2}$. $H_{S_2}$ is a disjoint union
  of less than $\k$ building blocks whose domains' inverse images decompose
  $\dom f$ into less than $\k$ disjoint pieces and one of them, say $C$, is not meager.
  By the Property of Baire one can find a basic open set $U$ such that $C\cap U$ is co-meager in $U$.
  Let $C(f)$ be a co-meager set in which $f$ is continuous.
  Now $f\rest (U\cap C\cap C(f))$ is a continuous reduction from $\B$ restricted to $(U\cap C\cap C(f))^2$
  to a building block of $H_{S_2}$. Thus it is sufficient to show 
  that this correctly chosen building block of $H_{S_1}$ is not reducible to any of the building blocks 
  of~$H_{S_2}$ on any such $U\cap C\cap C(f)$. 
  This will follow from Theorem~\ref{thm:nonRed} and Remark \ref{rem:nonRedbar}
  once we go through all the possible cases. So the following Lemma concludes the proof.

  \begin{Lemma}
    Assume that $\mu_1\in r$ is the least cardinal such that $(S_2\setminus S_1)\cap S^\k_{\mu_1}$
    is stationary. If $\mu_1>\o$, then
    \begin{myAlphanumerate}
    \item[(i)] for all $\mu_2>\o$, $\B_{\mu_1}(S_1)\not\le_B \B_{\mu_2}(S_2)$,
    \item[(ii)] $\B_{\mu_1}(S_1)\not\le_B \B^0_\o(S_2)$,
    \item[(iii)] $\B_{\mu_1}(S_1)\not\le_B \B^1_\o(S_2)$,
    \end{myAlphanumerate}
    and if $\mu_1=\o$, then
    \begin{myAlphanumerate}
    \item[(i*)] for all $\mu_2>\o$, $\B^0_{\o}(S_1)\not\le_B \B_{\mu_2}(S_2)$,
    \item[(ii*)] for all $\mu_2>\o$, $\B^1_{\o}(S_1)\not\le_B \B_{\mu_2}(S_2)$,
    \item[(iii*)] either 
      $$\B^0_{\o}(S_1)\not\le_B \B^0_\o(S_2)\text{ and }\B^0_{\o}(S_1)\not\le_B \B^1_\o(S_2)\eqno(1)$$
      or
      $$\B^1_{\o}(S_1)\not\le_B \B^0_\o(S_2)\text{ and }\B^1_{\o}(S_1)\not\le_B \B^1_\o(S_2).\eqno(2)$$
    \end{myAlphanumerate}
  \end{Lemma}
  \begin{proof}[Proof of the lemma] First we assume $\mu_1>\o$.
    \begin{myItemize}
    \item[(i)] There are two cases:
      \begin{myItemize}
      \item[Case 1:] $\mu_2=\mu_1$. Denote $B=A_{\mu_1}=A_{\mu_2}$ and
        $S_1'=(S_1\cap S^\k_{\mu_1})\cup B$ and $S_2'=(S_2\cap S^\k_{\mu_2})\cup B$.
        Now $\B_{\mu_1}(S_1)=\id\otimes E_{S_1'}$ and $\B_{\mu_2}(S_2)=\id\otimes E_{S_2'}$.
        Since by definition $B=B^\mu_{\nr}(K_\mu)$ where $K_\mu\subset S^\k_\o$ is stationary, and 
        $(S_2\setminus S_1)\cap S^\k_{\mu_1}$
        is stationary, the sets $S_1'$ and $S_2'$ satisfy the assumptions of 
        Theorem~\ref{thm:nonRed}.\ref{nr2b}, so the statement
        follows from Theorem~\ref{thm:nonRed}.\ref{nr2b} and Remark \ref{rem:nonRedbar}.
      \item[Case 2:] $\mu_2\ne \mu_1$. Let $S_1'=(S_1\cap S^\k_{\mu_1})\cup A_{\mu_1}$
        and $S_2'=(S_2\cap S^\k_{\mu_2})\cup A_{\mu_2}$ whence 
        $B_{\mu_1}(S_1)=\id\otimes E_{S_1'}$ and $B_{\mu_2}(S_2)=\id\otimes E_{S_2'}$.
        Now $S_1'\subset S^\k_{\ge\o}$ and $S_2'\subset S^\k_{\ge\o}$ 
        and since $A_{\mu_1}\cap A_{\mu_2}=\es$, the result 
        follows from Theorem \ref{thm:nonRed}.\ref{nr1b} and Remark \ref{rem:nonRedbar}.
      \end{myItemize}
    \item[(ii)] Let  $S_1'=(S_1\cap S^\k_{\mu_1})\cup A_{\mu_1}$, $S_2'=A^2_\o\cup ((S_2\cap S^\k_\o)\setminus A^0_\o)$,
      and $A_2'= A^0_\o$.
      By definition, 
      $$B^0_\o(S_2)=\id_{2^\k}\otimes E_{S_2'}\otimes E_{A_2'}$$
      and $B_{\mu_1}(S_1)=E_{S_1'}$. Since $A_{\mu_1}\cap A^2_\o=\es$, $S_1'\cap S^\k_\o=A_{\mu_1}$ 
      and $A^2_\o\subset S_2'$,
      we have that $S_2'\setminus S_1'$ is $\o$-stationary, because it contains $A^2_\o$.
      Also $A^0_\o\setminus S_1'=A^0_\o$, because $S_1'\cap A^0_\o=\es$,
      so $A_2'\setminus S_1'$ is $\o$-stationary. Now the result 
      follows from Theorem~\ref{thm:nonRed}.\ref{nr4b} and Remark \ref{rem:nonRedbar}.
    \item[(iii)] Similar to (ii).
    \end{myItemize}
    Then we assume $\mu_1=\o$.
    \begin{myItemize}
    \item[(i*)] Let $S_1'=A^2_\o\cup((S_1\cap S^\k_\o)\setminus A^0_\o)$, $A_1'=A^0_\o$ 
      $A_2'= A_{\mu_2}$ and $S_2'=(S_2\cap S^\k_{\mu_2})$. Since $A^0_\o\cap A_{\mu_2}=\es$, 
      we have that $A_2'\setminus A_1'$ is $\o$-stationary, so by 
      Theorem~\ref{thm:nonRed}.\ref{nr5} and Remark \ref{rem:nonRedbar}, 
      $$\id\otimes E_{S_1'}\otimes E_{A_1'}\not\le_B \id\otimes E_{S_2'\cup A_2'},$$
      which by definitions is exactly the subject of the proof.
    \item[(ii*)] Similar to (i*).
    \item[(iii*)] The situation is split into two cases, the latter of which is split into two subcases:
    \begin{myItemize}
    \item[Case 1:] \emph{$((S_2\setminus S_1)\cap S^\k_\o)\setminus (A^2_\o\cup A^0_\o)$ is stationary.}
      Let $S_1'=A^2_\o\cup((S_1\cap S^\k_\o)\setminus A^0_\o)$, $A_1'=A^0_\o$,
      $S_2'=A^2_\o\cup((S_2\cap S^\k_\o)\setminus A^0_\o)$ and $A_2'=A^0_\o$.
      Now $A'_2\setminus S_1'$ is obviously $\o$-stationary, since it is equal to $A^0_\o$.
      Also $S_2'\setminus S_1'$ is stationary, because
      it equals to $((S_2\setminus S_1)\cap S^\k_\o)\setminus (A^2_\o\cup A^0_\o)$ 
      which is stationary by the assumption.
      Now the first part of (1) follows from 
      Theorem~\ref{thm:nonRed}.\ref{nr3b} and Remark \ref{rem:nonRedbar},
      because 
      $\B^0_\o(S_1)=\id\otimes E_{S_1'}\otimes E_{A_1'}$ and 
      $\B^0_\o(S_2)=\id\otimes E_{S_2'}\otimes E_{A_2'}$.
      On the other hand let  
      $S_2''=A^3_\o\cup((S_2\cap S^\k_\o)\setminus A^1_\o)$ and $A_2''=A^1_\o$.
      Now $S_2''\setminus A_1'$ is stationary, because $A^3_\o\subset S_2''$ but $A^3_\o\cap A_1'=A^3_\o\cap A^0_\o=\es$.
      Also $A_2''\setminus A_1'$ is stationary since $A_2''\cap A_1'=A^1_\o\cap A^0_\o=\es$.
      Now also the second part of (1) follows from Theorem~\ref{thm:nonRed}.\ref{nr3b}  and Remark \ref{rem:nonRedbar},
      because 
      $B^0_1(S_1)=\id\otimes E_{S_1'}\otimes E_{A_1'}$ and $B^1_1(S_2)=\id\otimes E_{S_2''}\otimes E_{A_2''}$.      
    \item[Case 2:] \emph{$((S_2\setminus S_1)\cap S^\k_\o)\setminus (A^2_\o\cup A^0_\o)$ is non-stationary.}
      \begin{myItemize}
      \item[Case 2a:] \emph{$((S_2\setminus S_1)\cap S^\k_\o)\setminus (A^3_\o\cup A^1_\o)$ is stationary.}
        Now (2) follows from Theorem~\ref{thm:nonRed}.\ref{nr3b} and 
        Remark \ref{rem:nonRedbar} in a similar way as (1) 
        followed in Case~1.
      \item[Case 2b:] \emph{$((S_2\setminus S_1)\cap S^\k_\o)\setminus (A^3_\o\cup A^1_\o)$ is non-stationary.}
        Now we have both:
        $$((S_2\setminus S_1)\cap S^\k_\o)\setminus (A^2_\o\cup A^0_\o)\text{ is non-stationary}\eqno(*)$$
        and
        $$((S_2\setminus S_1)\cap S^\k_\o)\setminus (A^3_\o\cup A^1_\o)\text{ is non-stationary.}\eqno(**)$$
        Now from $(*)$ it follows that $S_2\setminus S_1\subset_{\NS(\o)}A^2_\o\cup A^0_\o$.
        From $(**)$ it follows that $S_2\setminus S_1\subset_{\NS(\o)} A^3_\o\cup A^1_\o$.
        This is a contradiction, because $S_2\setminus S_1$ is $\o$-stationary and 
        $(A^2_\o\cup A^0_\o)\cap (A^3_\o\cup A^1_\o)=\es$.
      \end{myItemize}
    \end{myItemize}
  \end{myItemize}  
\end{proof} 
\end{proof}

\section{On Chains In $\la\E^B_\k,\le_B\ra$}

There are chains of order type $\k^+$ in Borel equivalence relation on $2^\k$:

\begin{Thm}\label{thm:Chain} 
  Let $\k>\o$.
  There are equivalence relations $R_i\in \E^B_\k$, for $i<\k^+$, such that
  $i<j\iff R_i\lneq_B R_j\lneq E_0$.
\end{Thm}
\begin{RemarkN}\label{rem:saturated}
  In many cases there are $\k^+$-long chains in the power set of $\k$ ordered by inclusion 
  modulo the non-stationary ideal whence a weak version of 
  this theorem could be proved using Theorem~\ref{thm:Main2}. 
  Namely if the ideal $I^\k_{\NS}$ of non-stationary subsets of $\k$
  is \emph{not} $\k^+$-saturated, then there are $\k^+$-long chains. 
  In this case being \emph{not} $\k^+$-saturation means that 
  there exists a sequence $\la A_i\mid i<\k^+\ra$ of subsets of $\k$
  such that $A_i$ is stationary for all $i$ but $A_i\cap A_j$ is non-stationary for all $i\ne j$.
  Now let $f_\a$ be a bijection from $\k$ to $\a$ for all $\a<\k^+$ and let 
  $$B_\a=\mathop{\triangledown}_{i<\a}A_i=\{\a\mid\text{ for some }i<\a, \a\in A_{f_\a(i)}\}$$
  It is not difficult to see that $\la B_\a \mid \a<\k^+\ra$ is a chain.
  On the other hand the existence of such a chain implies that $I^{\k}_{\NS}$ is not $\k^+$-saturated.
  
  By a theorem of Gitik and Shelah \cite[Theorem 23.17]{4Jech}, $I^\k_{\NS}$ is not $\k^+$-saturated
  for all $\k\ge\aleph_2$. 
  By a result of Shelah \cite[Theorem 38.1]{4Jech}, it is consistent relative to 
  the consistency of a Woodin cardinal that $I^{\aleph_1}_{\NS}$ \emph{is}
  $\aleph_2$-saturated in which case there are no chains of length $\o_2$ in
  $\la \Po(\o_1),\subset_\NS\ra$. 
  On the other hand in the model provided by Shelah, CH fails.
  According to Jech \cite{5HandBookOfSetTheory} it is an open question whether CH implies
  that $I^{\aleph_1}_{\NS}$ is not $\aleph_2$-saturated.
  
  However, as the following shows, it follows from ZFC that 
  there are $\k^+$-long chains in $\la\E^B_\k,\le_B\ra$ for any uncountable $\k$.
\end{RemarkN}
\begin{proof}[Proof of Theorem \ref{thm:Chain}]
  By the proof of Corollary \ref{cor:Antichain}, page~\pageref{cor:Antichain}, one can find $\o$-stationary
  sets $S_i$ for $i<\k^+$ such that $S_i\setminus S_j$ and
  $S_j\setminus S_i$ are stationary whenever $i\ne j$. For all $j\in [1,\k^+)$, let 
  $$R_j=\bigoplus_{i<j} E_{S_i},$$
  where the operation $\oplus$ is from Definition~\ref{def:Plus}, page~\pageref{def:Plus}. 

  Let us denote $P_{A}=\Cup_{i\in A}2^{\k\times\{i\}}$ for $A\subset \k^+$, i.e. for example
  $P_{j}=\Cup_{i<j}2^{\k\times\{i\}}$.

  Let us show that
  \begin{myEnumerate}
  \item\label{chain1} if $i<j$, then $R_i\le_B R_j$,
  \item\label{chain2} if $i<j$, then $R_j\not\le_B R_i$,
  \item\label{chain3} for all $i<\k^+$, $R_i\lneq_B E_0$.
  \end{myEnumerate}
  Item \ref{chain1} is simple: let $f\colon P_i\to P_j$ be the inclusion map (as $P_i\subset P_j$).
  Then $f$ is clearly a reduction from $R_i$ to $R_j$. 
  
  Suppose then that $i<j$ and that $i\le k < j$. To prove \ref{chain2} it is sufficient to show that there 
  is no reduction from $E_{S_k}$ to $R_j$.  
  Let us assume that $f\colon 2^\k\to P_j$ is a Borel reduction from 
  $E_{S_k}$ to $R_j$.
  Now 
  $$2^\k=\Cup_{\a<i} f^{-1}[P_{\{\a\}}],$$
  so one of the sets $f^{-1}[P_{\{\a\}}]$ is not meager; let $\a_0$ be an index witnessing this.
  Note that $\a_0<k$, because $\a_0<i\le k$.
  Because $f$ is a Borel function and Borel sets have the Property of Baire,
  we can find a $p\in 2^{<\k}$ such that $C=N_p\cap C(f)\cap f^{-1}[P_{\{j\}}]$ is co-meager in~$N_p$.
  But now $f\rest C$ is a continuous reduction from $E_{S_k}\cap C^2$ to $E_{S_\a}$
  which contradicts Theorem~\ref{thm:nonRed}.\ref{nr1b}.

  To prove \ref{chain3} we will show first that $R_i\le_B \bigoplus_{j<i}E_0$ and then
  that $\bigoplus_{j<i}E_0\le_B E_0$, after which we will show that $E_0\not\le_B R_i$ for all $i$. 

  Let $f_j$ be a reduction from $E_{S_j}$ to $E_0$ for all $j<i$ given by Claim 3 of
  the proof of Theorem \ref{thm:Main1}. Then combine these reductions to get a reduction 
  from $R_i$ to $\bigoplus_{j<i}E_0$. To be more precise,
  for each $\eta\in P_{\{j\}}$ let $f(\eta)$ be $\xi$ such that $\xi\in P_{\{j\}}$ and 
  $\xi=f_{j}(\eta)$.

  Let $\{A_k\mid k\le i\}$ be a partition of $\k$ into disjoint unbounded sets.
  Let $\eta\in P_{i}$. By definition, $\eta\in P_{\{k\}}$ for some $k<i$.
  Define $\xi=F(\eta)$ as follows. Let $f\colon A_i\to \k$ be a bijection.
  \begin{myItemize}
  \item If $\a\in A_i$, then let $\xi(\a)=\eta(f(\a))$.
  \item If $\a\in A_j$ and $j\ne k$, then let $\xi(\a)=0$.
  \item If $\a\in A_k$, then let $\xi(\a)=1$.
  \end{myItemize}
  It is easy to see that $F$ is a continuous reduction.  

  Assume for a contradiction that $E_0\le_B R_i$ for some $i<\k^+$. Then by \ref{chain1} and
  transitivity, $E_0\le_B R_j$ for all $j\in [i,\k^+)$. By the above also $R_j\le_B E_0$ for all $j\in [i,\k^+)$
  which, again by transitivity, implies that the relations $R_j$ for $j\in [i,\k^+)$ are mutually bireducible
  to each other which contradicts~\ref{chain2}.
\end{proof}

\newpage

\label{ch:thebibl}

\end{document}